\documentclass[11pt]{article}

\usepackage{amsmath}
\usepackage{verbatim}
\usepackage{amssymb}
\usepackage{latexsym}
\usepackage{amsthm,dsfont}
\usepackage[T1]{fontenc}
\usepackage[utf8]{inputenc}

\usepackage{mathrsfs}
\usepackage{delarray}
\usepackage{pstricks}
\usepackage{pst-plot}
\usepackage{pst-all}
\usepackage{auto-pst-pdf}
\usepackage{biblatex}
\usepackage{graphicx}
\usepackage{hyperref}

\bibliography{bib_nudging}

\usepackage{enumerate}
\newtheorem{Theorem}{Theorem}

\newtheorem{Definition}[Theorem]{Definition}

\newtheorem{Remark}[Theorem]{Remark}

\def\eqn{\begin{equation}}
\def\nqe{\end{equation}}

\def\eqa{\begin{eqnarray}}
\def\aqe{\end{eqnarray}}

\def\Ne{\mathds{N}}

\def\Re{\mathds{R}}

\def\1e{\mathds{1}}

\def\vg0{\mathbf{0}}

\def\C{\mathscr{C}}
\def\D{\mathscr{D}}

\def\G{\mathscr{G}}
\def\H{\mathscr{H}}
\def\I{\mathscr{I}}

\def\L{\mathscr{L}}
\def\M{\mathscr{M}}

\def\O{\mathscr{O}}
\def\P{\mathscr{P}}

\def\S{\mathscr{S}}
\def\T{\mathscr{T}}
\def\U{\mathscr{U}}
\def\V{\mathscr{V}}
\def\W{\mathscr{W}}
\def\X{\mathscr{X}}

\def\d{\displaystyle}

\def\e{\varepsilon}

\def\scal#1#2{\langle #1,#2 \rangle}

\def\set#1#2{\left\{\mskip 1mu #1 \mskip 1mu
    | \mskip 1mu #2 \mskip 1mu \right\}
    }

\def\mod#1{|\mskip 1mu #1 \mskip 1mu|}

\def\norm#1{\left\| \mskip 1mu #1 \mskip 1mu \right\|}

\def\operator#1#2#3#4#5{
    \begin{array}{lcll}
    \d #1 \colon & \d #2 & \longrightarrow & \d #3 \\[.4ex]
                 & \d #4 & \longmapsto     & \d #5
    \end{array}
    }

\def\twm1{T_{\!W}^{-1}}

\def\ud{\,\mathrm{d}}

\begin{document}

\parindent0pt
\parskip8pt
%-------------------
\title{{\bf Application of a nudging technique to thermoacoustic
tomography}}
\author{Xavier Bonnefond and S\'ebastien Marinesque}

%\date{}
%-------------------
\maketitle

%%%%%%%%%%%%%%%%%%%%%%%%%%%%%%%%%%%%%%%%%%%%%%%%%%%%%%%%
%%%%%%%%%%%%%%%%%%%%%%%%%%%%%%%%%%%%%%%%%%%%%%%%%%%%%%%%
\begin{abstract}
ThermoAcoustic Tomography (TAT) is a promising, non invasive, medical imaging technique whose inverse problem can be formulated as an initial condition reconstruction. In this paper, we introduce a new algorithm originally designed to correct the state of an evolution model, the \emph{back and forth nudging} (BFN), for the TAT inverse problem. We show that the flexibility of this algorithm enables to consider a quite general framework for TAT. The backward nudging algorithm is studied and a proof of the geometrical convergence rate of the BFN is given. A method based on Conjugate Gradient (CG) is also introduced. Finally, numerical experiments validate the theoretical results with a better BFN convergence rate for more realistic setups and a comparison is established between BFN, CG and a usual inversion method.
\end{abstract}

%%%%%%%%%%%%%%%%%%%%%%%%%%%%%%%%%%%%%%%%%%%%%%%%%%%%%%%%%%%%%
\section{Introduction}%%%%%%%%%%%%%%%%%%%%%%%%%%%%%%%%%%%%%%
%%%%%%%%%%%%%%%%%%%%%%%%%%%%%%%%%%%%%%%%%%%%%%%%%%%%%%%%%%%%%%
\label{introduction}

ThermoAcoustic Tomography (TAT) is a hybrid imaging technique that uses ultrasound waves produced by a body submitted to a radiofrequency pulse, uniformly deposited
throughout the body. The absorption of this initial energy causes a non-uniform thermal expansion, leading to the propagation of a pressure wave
outside the body to investigate. This wave is then measured all around the body with piezoelectric transducers or, more recently, thanks to interferometry 
techniques. 

It appears that the absorption of the initial pulse is highly related to the physiological properties of the tissue~\cite{joines}. As a result, the magnitude of 
the ultrasonic emission ({\sl i.e.} thermoacoustic signal), which is proportional to the local energy deposition, reveals physiologically specific absorption contrast. See
\cite{kruger, kruger2, paltauf_interfero, wang_photoacoustic} for an introduction to the experimental setup. Considering that the initial illumination is a Dirac distribution in time, the problem of recovering the absorptivity of the investigated body from the thermoacoustic signal is equivalent to recovering the initial condition of a Cauchy problem involving the wave equation from the knowledge of the solution on a surface surrounding the imaging object~\cite{otmar_book}. Note that the experimental constraints do not allow the acquisition surface to completely surround the tissues to investigate, so that one cannot expect a better situation than measurements on a half-sphere (as in breast cancer detection for example). Moreover, several experimentation constraints limit the number of sensors, so that the methods of reconstruction considered must be robust to the space sampling.

Many works dealing with TAT have been achieved during the last decade (see~\cite{kuchment} for an overview, or~\cite{akk,wangbook}), and many of the methods used in these papers strongly depend on a set of additional assumptions, among which:
\begin{itemize}
\item the \emph{homogeneity of the tissues}, leading to a constant speed of sound for the pressure wave;
\item the \emph{lack of frequency-dependent attenuation}, so that the pressure wave obeys the classical, undamped,
wave equation;
\item the \emph{complete data situation}, where the acquisition surface is considered to enclose the whole body.
\end{itemize}

The study of TAT in the case of an homogeneous medium without attenuation gave rise to explicit inversion formulas. These latter usually require that the acquisition set is a closed surface, even if they can be approximated in the \emph{limited view} situation (where the surface is a half-sphere for example). Among these, filtered backprojection formulas, even though they have been very successfully implemented in several situations, seem to face some issues when the source has support partly outside the observation surface, or when the data sampling is not sufficient (see~\cite{finch, otmar_fbp}). Fourier's type formulas, which offer a much better numerical efficiency, apply on very specific geometries (see~\cite{kunyansky_serie,norton_3d, wang_timedomain2}) or imply some additional approximations, like an interpolation from spherical to cartesian coordinates~\cite{anastasio_fourier, cox_cavite, otmar_bessel, haltmeier_fftnonuni, wang_timedomain1, kuchment_limited}.
Even though some of these techniques give rise to very efficient reconstruction schemes (especially in~\cite{kunyansky_serie}), so far it is not clear if they can be extended to less restrictive acquisition geometries or to an attenuated wave with non constant speed (see~\cite{kunyansky_curve} for an answer about the acquisition surface).

Moreover, since the work of Burgholzer \textit{et al} (see~\cite{paltauf_time_reversal}), the \emph{time reversal} method (first suggested in~\cite{finch}) has been applied with success to the TAT (see \textit{e.g.}~\cite{hristova,kuchment_time_reversal}). Even though the first known results on this method required a complete data framework, recent works by Stefanov \textit{et al} extended it to a Neumann Series method in a fairly general situation (incomplete data, known variable sound speed and external source) with very good results in~\cite{stefanov_algorithm,stefanov}. 

%% Thanks to the relation between solutions of the classical wave equation and spherical means~\cite{hilbert},
%% this situation (constant speed, no attenuation) allows the use of variational methods. 
%% A regularization scheme is proposed in~\cite{xapi},
%% based on the \emph{approximate inverse} (see~\cite{schuster} and~\cite{otmar_fbp} for direct applications of this concept), which
%% remains stable in the limited view situation. Yet, any attempt to weaken the assumptions would certainly 
%% face serious technical difficulties.

In this paper, we show that the framework of TAT lends itself to the application of some \emph{data assimilation} techniques. These latter, mainly used in geophysics, aim at correcting the state of an evolutionary model by means of data, in order to obtain a good approximation of the real state (see~\cite{Tal1} for an introduction). The particular data assimilation method proposed here is based on a \emph{nudging} technique: given an evolution model of the state and direct observation (our data), it consists in adding, inside the model equation, a newtonian recall of the state solution to the observations (or data), which is usually called the feedback or nudging term~\cite{AurBl3}. As we shall see, from a practical point of view, this method can be successfully used to manage the usual issues of the TAT inverse problem as incomplete data, external source and variable sound speed (when given, however). So far, however, the theoretical convergence result for the nudging technique is based on a classical result about stabilization of the wave equation~\cite{lions,russel}, which requires somehow a geometric optics condition (see~\cite{bardos} or~\cite{liu97} for an overview). This geometric condition, whose verification is not a simple matter, is not automatically satisfied in an incomplete data framework or with a variable sound speed. Consequently the proof provided in this article only stands in a very favorable situation (3d case, constant speed and complete data), but does not depend on usual stabilization techniques and yields explicit convergence rates.

Different nudging terms can be found in the literature, and most of them have been used to assimilate data in physical oceanography as in~\cite{Krishnamurti,Lorenc,Lyne} and more recently in~\cite{AurouxetBonnabel}. Fundamental articles offer the basis of most popular techniques as Kalman filter~\cite{Kalman}, Kalman-Bucy filter~\cite{KB} and Luenberger observer~\cite{Luen}.

Since the basic nudging method showed its limits in real conditions by using the data once, Auroux and Blum proposed to extend the method by adding a resolution backward in time and iterating the process, thus defining the \emph{back and forth nudging} (BFN) algorithm. They presented this technique in~\cite{AurBl1,AurBl3} and implemented it with computationally efficient results in comparison with traditional data assimilation methods, like optimal filters (from extended Kalman's to particle filters) or other optimal minimization techniques (variational methods), which appear costfull respectively in memory needs and computing time -- see the brief introductions of Talagrand in~\cite{Tal1,Tal2,Tal3} for further details about data assimilation and a relevant bibliography. Recently, a more general formulation of the BFN using observers has been presented by Ramdani, Tucsnak and Weiss in~\cite{ramdani}.

As inverse problems are often solved by means of variational techniques, we also introduce a least squares method based on the Conjugate Gradient (CG).

The two methods introduced are compared to the Neumann Series presented in~\cite{stefanov_algorithm}, which is one of the best existing method which fits the variable speed case and partial data settings.

This paper is organized as follows. In Section~\ref{pres}, we introduce the TAT inverse problem in a quite general form, {\sl i.e.} without the use of the usual assumptions, then we describe the BFN algorithm and we state the main result of this article: the convergence Theorem~\ref{theoreme_principal}. Then Section~\ref{section_convergence} is devoted to the proof of this result. In Section~\ref{var}, we describe the variational method that uses CG algorithm. Finally, we give numerical results in Section~\ref{section_illustration}.

\subsection*{Notation}

In the following, we shall use this notation:
\begin{itemize}
\item The open ball with center~$x\in\Re^3$ and radius~$r>0$ is denoted by~$B(x,r)$;
\item Its boundary, the sphere with center~$x$ and radius~$r$, is noted~$S(x,r)$;
\item For every~$x\in\Re^3$ and~$0<r_1<r_2$, we define the spherical shell of center~$x$ and radii~$r_1,\ r_2$:
$$
A(x,r_1,r_2)=\set{y\in\Re^3}{r_1<\norm{x-y}<r_2},
$$
and if~$x=0$:
$$
A(r_1,r_2)=A(0,r_1,r_2);
$$
\item For any set~$S$ of~$\Re^3$ and~$r_1>0$, we denote by~$T(S,r_1)$ the set~$S$ extended to the points at a distance lower than~$r_1$ from~$S$:
$$
T(S,r_1)=\set{y\in\Re^3}{\exists x \in S, \norm{y-x}<r_1}.
$$
\end{itemize}

%%%%%%%%%%%%%%%%%%%%%%%%%%%%%%%%%%%%%%%%%%%%%%%%%%%%%%%%%%%%%%%%%%%%%%%%%%%%
\section{Presentation of the method}
\label{pres}%%%%%%%%%%%%%%%%%%%%%%%%%%%%%%%%%%%%%%
%%%%%%%%%%%%%%%%%%%%%%%%%%%%%%%%%%%%%%%%%%%%%%%%%%%%%%%%%%%%%%%%%%%%%%%%%%%%
\subsection{The general TAT problem}
\label{general_problem}

In the following, we will denote by~$p_{\mathrm{exact}}(x,t)$ the pressure wave resulting from the thermal expansion of the body. We will make the assumption that the measurement process of this pressure is subject to some perturbation~$p_ \mathrm{noise}$, such that the actual data can be written~$p_\mathrm{data}:=p_\mathrm{exact}+p_\mathrm{noise}$. In thermoacoustic tomography, the data are nothing but~$\set{p_{\mathrm{data}}(x,t)}{x\in \S}$, where~$\S$ is a surface surrounding the body to investigate. The pressure~$p_\mathrm{exact}$ satisfies the Cauchy problem (see \cite{otmar_book} for a detailed calculation):
$$
\begin{array}|{ll}.
L p_{\mathrm{exact}}(x,t)=E_{dep}(x,t), & (x,t)\in\Re^3\times[-1,\infty),\\
p_{\mathrm{exact}}(x,-1)=0, & x\in\Re^3,\\
\partial_t p_{\mathrm{exact}}(x,-1)=0,& x\in\Re^3,
\end{array}
$$
where~$L$ is the differential operator governing the wave, most likely the wave operator, or a damped wave
operator, and~$E_{dep}$ is the energy deposited in the body around the time~$t=0$ (so we have to consider a negative initial time, {\sl e.g.}~$t=-1$). This energy can be written:
$$
E_{dep}=e(x)\dfrac{\ud j}{\ud t}(t),
$$
$$
e(x):=\dfrac{\beta c}{c_p}I_{em}(x)\Psi(x),
$$
where~$e$ is called the \textit{normalized energy deposition function}. Here~$\beta$ is the
\textit{thermal expansion coefficient},~$c_p$ is the \textit{specific heat capacity},~$c$ the sound speed
and~$I_{em}$ is the \textit{radiation intensity} of the initial energy pulse. All these parameters, 
including the time shape~$j$ of the pulse, are known and we see that the knowledge of~$E_{dep}$ is 
sufficient to compute the absorption density~$\Psi$ inside the body, which is our purpose.

So far, the general thermoacoustic problem can be formulated in the following way:

\begin{quotation}
\sf
Let~$p_{\mathrm{exact}}(x,t)$ be solution of:
\eqn
\label{eq_onde_generale}
\begin{array}|{ll}.
L p_{\mathrm{exact}}(x,t)=f_0(x)\dfrac{\ud j}{\ud t}(t), & (x,t)\in\Re^3\times[-1,\infty), \\
p_{\mathrm{exact}}(x,-1)=0, & x\in\Re^3, \\
\partial_t p_{\mathrm{exact}}(x,-1)=0,& x\in\Re^3,
\end{array}
\nqe

where~$j$ is known and~$f_0$ is the object to reconstruct, with compact support in the unit ball~$B(0,1)$. Can we compute~$f_0$, or a good
approximation of~$f_0$, from the knowledge of~$p_{\mathrm{data}}=p_\mathrm{exact}+p_\mathrm{noise}$ on a surface~$\S$ surrounding~$B(0,1)$ ? 
\end{quotation}

%%%%%%%%%%%%%%%%%%%%%%%%%%%%%%%%%%%%%%%%%%%%%%%%%%%%%%%%%%%%%%%%%%%%%%
\subsection{Reduction to the homogeneous case}

As we have seen, even though it is not possible in practice, the initial energy deposition is meant
to be a Dirac pulse in time at~$t=0$, that is~$\frac{\ud j}{\ud t}=\delta'_0$. Let~$q_0$ be solution, which is no longer to be considered for negative time values, where it is the null solution, of the Cauchy problem:
$$
\begin{array}|{ll}.
L q_0(x,t)=f_0(x)\delta'_0(t), & (x,t)\in\Re^3\times\Re_+,\\
q_0(x,0)=0,  & x\in\Re^3, \\
\partial_t q_0 (x,0)=0, & x\in\Re^3,
\end{array}
$$
which is equivalent to:
\eqn
\label{eq_onde}
\begin{array}|{ll}.
L q_0(x,t)=0, & (x,t)\in\Re^3\times\Re_+, \\
q_0(x,0)=f_0(x),& x\in\Re^3,\\
\partial_t q_0(x,0)=0,& x\in\Re^3,
\end{array}
\nqe
according to Duhamel's principle (see \cite{hilbert}), provided that~$L$ is a (damped) wave operator. Then~$p_{\mathrm{data}}$ is solution
of \eqref{eq_onde_generale} if and only if~$p_{\mathrm{data}}=j *_t q_0$, which means that a deconvolution operation
on the thermoacoustic signal~$p_{\mathrm{exact}}$ leads to the knowledge of the solution~$q_0$ of \eqref{eq_onde} on the
same surface~$\S$ surrounding the object.

In the following, we will still denote by~$p_{\mathrm{exact}}$ the solution of \eqref{eq_onde}, assuming that we know this
latter for all~$x$ in~$\S$, and all~$t$ in~$\Re_+$.

%%%%%%%%%%%%%%%%%%%%%%%%%%%%%%%%%%%%%%%%%%%%%%%%%%%%%%%%%%%%%%%%%%%%%%%%%%%%%%
\subsection{Some useful facts about the standard wave equation}
 \label{useful_facts}
In this section we give some basic properties of the following wave equation:
\eqn
\label{eq_onde_exemple}
\begin{array}|{ll}.
L_0 u(x,t)=0, & (x,t)\in\Re^3\times\Re_+,\\
u(x,0)=l(x), & x\in\Re^3,\\
\partial_t u(x,0)=h(x),& x\in\Re^3,
\end{array}
\nqe
where~$l$ and~$h$ are two~$C^{\infty}$ functions with compact support in~$B(0,1)$, and~$L_0:=\partial_{tt}-\Delta$. All this material can be found in~\cite{john} for example. According to Huyghens' principle in the strong sense, for every couple~$(x_0,t_0)\in\Re^3\times\Re_+$, the classical solution~$u(x_0,t_0)$ of~\eqref{eq_onde_exemple} only depends on the values of~$l$ and~$h$ on the
backward characteristic cone:
$$
\norm{x-x_0}=\mod{t-t_0}.
$$

As a matter of fact, when the initial conditions~$l$ and~$h$ have their support in~$B(0,1)$, the solution at any time~$\tilde t$ has its support in~$B(0,1+\tilde t)$.

Huyghens' principle leads to the energy conservation for the solution~$u$ of~\eqref{eq_onde_exemple}.
 
\begin{Definition}\sf
For every time~$T$, we call \textit{energy} of~$u$ the quantity:
$$
E_u(T):=\dfrac{1}{2}\left(\norm{u(\cdot,T)}_{H^1_0(\Re^3)}^2+
\norm{\partial_t u(\cdot,T)}_{L^2(\Re^3)}^2 \right ).
$$
\end{Definition}

The variation of the energy of~$u$ is estimated from the following identity:
$$
\partial_t u \left (\partial_{tt} u-\Delta u \right )=0=\partial_t\left [ \dfrac{\norm{\nabla u}^2}{2}+\dfrac{(\partial_t u)^2}{2} \right ]
-\sum_{i=1}^{3}\partial_{x_i}(\partial_t u\, \partial_{x_i} u).
$$
Since~$u$ has a compact support, integrating for~$0<t<T$ and over the whole space yields:
$$
E_u(T)=E_u(0).
$$
Note that if~$h$ is zero, we have:
$$
E_u(0)=\dfrac{1}{2}\norm{l}_{H^1_0(\Re^3)}^2.
$$

%%%%%%%%%%%%%%%%%%%%%%%%%%%%%%%%%%%%%%%%%%%%%%%%%%%%%%%%%%%%%%%%%%%%%%%%%%%%
\subsection{Nudging}
\label{nudging}

Now, let us introduce the application of the BFN algorithm to the general TAT problem.

\begin{figure}[ht]
\centering
% Generated with LaTeXDraw 2.0.8
% Thu Sep 30 16:37:11 CEST 2010
% \usepackage[usenames,dvipsnames]{pstricks}
% \usepackage{epsfig}
% \usepackage{pst-grad} % For gradients
% \usepackage{pst-plot} % For axes
\scalebox{.8} % Change this value to rescale the drawing.
{
\begin{pspicture}(0,-4.249219)(9.061875,4.289219)
\definecolor{color0b}{rgb}{0.6,0.6,1.0}
\definecolor{color15b}{rgb}{1.0,0.4,0.4}
\pscircle[linewidth=0.02,dimen=outer,fillstyle=solid,fillcolor=color15b](4.05,-0.15921874){4.09}
\pscircle[linewidth=0.02,dimen=outer,fillstyle=solid](4.05,-0.15921874){3.51}
\psbezier[linewidth=0.04,fillstyle=solid,fillcolor=color0b](2.4395187,1.5323794)(3.0971763,2.5487037)(2.476185,1.8549681)(3.6527324,2.2281444)(4.8292794,2.6013205)(4.4320264,3.1307812)(5.410654,2.3960876)(6.3892817,1.661394)(7.24,-1.659478)(6.1534376,-2.2343483)(5.0668755,-2.8092186)(6.05046,-1.3635328)(4.8411865,-1.1067293)(3.6319132,-0.84992564)(3.1048105,-2.6530087)(2.1424053,-1.8984618)(1.18,-1.1439148)(1.7818613,0.516055)(2.4395187,1.5323794)
\usefont{T1}{ptm}{m}{n}
\rput(3.1059375,-0.51921874){Body}
\pscircle[linewidth=0.06,linestyle=dashed,dash=0.16cm 0.16cm,dimen=outer](4.05,-0.15921875){3.81}
\psline[linewidth=0.024](7.14,-2.3092186)(7.68,-3.0092187)(8.18,-3.0092187)
\usefont{T1}{ptm}{m}{n}
\rput(8.491406,-3){$\S$}
\psline[linewidth=0.024](5.6,3.5107813)(6.3,3.9907813)(6.68,3.9907813)
\usefont{T1}{ptm}{m}{n}
\rput(7.414063,3.95007814){$T(\S,\e)$}
\psbezier[linewidth=0.048,linestyle=dotted,dotsep=0.16cm,fillcolor=color0b](2.3,1.6907812)(3.0798693,2.8671925)(2.314355,1.9686368)(3.7,2.4507813)(5.085645,2.9329257)(4.459397,3.3707812)(5.52,2.4907813)(6.580603,1.6107812)(7.54,-1.6594647)(6.26,-2.3092186)(4.98,-2.9589727)(5.76,-2.1092188)(4.74,-1.5092187)(3.72,-0.9092187)(3.1287885,-2.9803994)(1.98,-1.9892187)(0.8312115,-0.99803823)(1.5201306,0.5143699)(2.3,1.6907812)
\psbezier[linewidth=0.036,linestyle=dotted,dotsep=0.16cm,fillcolor=color0b](1.82,1.6707813)(2.6598694,2.8471925)(2.4143548,2.4286368)(3.72,2.8907812)(5.0256453,3.3529258)(4.999397,3.3307812)(6.04,2.4907813)(7.080603,1.6507813)(7.78,-1.9594648)(6.5,-2.6692188)(5.22,-3.3789728)(5.88,-3.2692187)(4.92,-2.2892187)(3.96,-1.3092188)(2.9487884,-3.4003992)(1.74,-2.3892188)(0.5312115,-1.3780382)(0.9801307,0.4943699)(1.82,1.6707813)
\end{pspicture} 
}
\caption{\label{schema} The body to investigate is surrounded by the acquisition surface~$\S$, and generates a pressure wave.}
\end{figure}
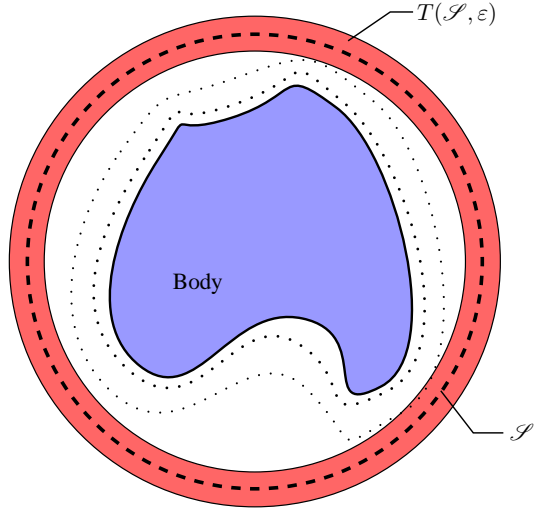

Assume that the observation surface~$\S$ is included in the sphere $S(0,1)$. As exposed in Section~\ref{general_problem}, our purpose is to compute an approximation of the
original object~$f_0$ from an incomplete set of data~$\phi p_{\mathrm{data}}$, where~$p_{\mathrm{data}}$ is a the solution (possibly contaminated by noise) of the wave equation \eqref{eq_onde} and where $\phi$ allows the knowledge of~$p_{\mathrm{data}}$ only on the observation surface $\S$. Theoretically speaking,~$\phi$ should be~$\1e_{\S}$, but in practice the discretization process allows us to choose a~$C^{\infty}$ function with compact support. Indeed, given a resolution, there exists~$\varepsilon >0$ such that any function supported in~$T(\S,\varepsilon)$ and equal to $1$ on~$\S$ has the same discrete counterpart as $\1e_{\S}$ (see Figure~\ref{schema}). This is the reason why we will make the following assumptions on~$\phi$:

\begin{enumerate}[i]
\item $\forall x \in \Re^3,\quad 0\leq\phi(x)\leq 1$;
\item $\phi$ has a compact support in~$T(\S,\varepsilon)$;
\item $\forall x \in T(\S,\varepsilon/2),\quad \phi(x)=1$.
\end{enumerate}

Finally, in order to avoid some technical difficulties, we will assume that the initial
object~$f_0$ has his support in~$B(0,1-\varepsilon)$ and is such that every Cauchy 
problem we shall consider has a classical solution, say~$C^\infty$.

We may proceed as follows:
\begin{Definition}[Iterate of the BFN algorithm for TAT]{\quad}
\label{def_bfn}
\vspace{-.5cm}
\begin{enumerate}
\item \,{\sf [Forward evolution]}
 
From an rough estimate~$f_{0,1}$ of the object to reconstruct~$f_0$, defined as an approximate of the initial condition supported in~$B(0,1-\varepsilon)$, we first consider a solution~$p_1$ of Problem~\eqref{eq_onde_exemple} with initial conditions~$l=f_{0,1}$ and~$h=0$ (as usual in TAT).

We compute this solution until time $T=2$, chosen so that both~$p_1(\cdot,t)$ and~$p_{\mathrm{exact}}(\cdot,t)$ (just as every function with compact support in the ball~$B(0,1-\varepsilon)$) vanish on~$B(0,1+\varepsilon)$ for all~$t\geq T$ (\textit{cf} Huyghens' principle in Section~\ref{useful_facts}).

Thus at this final time, in the noise-free situation, the \textit{innovation term}~$p_1-p_{\mathrm{data}}$, defined as the difference between the solution~$p_1$ of Problem~\eqref{eq_onde_exemple} and the data, is still a solution of~\eqref{eq_onde_exemple} with~$l=f_0-f_{0,1}$ and~$h=0$ and has an energy matching up to:
$$
E_{p_1-p_{\mathrm{data}}}(T)=E_{p_1-p_{\mathrm{data}}}(0)=\dfrac{1}{2}\norm{\nabla(f_0-f_{0,1})}_{H^1_0(\Re^3)}^2,
$$
since the conservation of energy for the wave equation holds in this situation.

\item \,{\sf [Backward evolution]}

Then we apply the backward nudging method to make the solution evolve back in time, with, as announced, the addition in the backward equation of a newtonian feedback, which adjusts the solution along its evolution by recalling it to the observed data.

Namely, we add to~$L$ the feedback correction term~$k \phi \partial_t(\cdot-p_\mathrm{data})$, for some nudging parameter~$k\in\Re_+$. And we compute the backward solution starting from the final state of the forward implementation, {\sl i.e.}~$p_1(x,T)$ and~$\partial_t p_1(x,T)$, untill the initial time~$t=0$.

We obtain a corrected solution~$\tilde p_1$ of the Cauchy problem, called \textit{backward nudging problem}:
\eqn
\label{eq_nudging_back0}
\begin{array}|{ll}.
L \tilde p_1(x,t)=k \phi(x) \partial_t (\tilde p_1-p_{\mathrm{data}})(x,t), & (x,t)\in\Re^3\times[0,T], \\
\tilde p_1(x,T)=p_1(x,T), & x\in\Re^3,\\
\partial_t \tilde p_1(x,T)=\partial_t p_1(x,T), & x\in\Re^3.
\end{array}
\nqe
In order to obtain an equation in forward time, let us consider the map~$t\mapsto \tilde p_1(T-t)$, still denoted by~$\tilde p_1$. The Cauchy problem becomes:
\eqn
\label{eq_nudging_back1}
\begin{array}|{ll}.
\tilde L \tilde p_1=-k \phi(x) \partial_t (\tilde p_1+p_{\mathrm{data}})
(x,T-t),& (x,t)\in\Re^3\times[0,T], \\
\tilde p_1(x,0)=p_1(x,T), & x\in\Re^3,\\
\partial_t \tilde p_1(x,0)=-\partial_t p_1(x,T),& x\in\Re^3,
\end{array}
\nqe
where~$\tilde L$ is the operator~$L$ backwards in time.

Obviously, in practice, when~$L$ is not the wave operator, there is no reason for this operator to be reversible. Nevertheless the newtonian feedback may act as a regularization term and keep the computation stable.

\item \,{\sf [Update of the estimate]}

Finally, after this back and forth evolution, a new estimate~$\tilde p_1(T)$ is obtained, but in order to iterate the process we need the new estimate to be supported in~$B(0,1-\varepsilon)$, just like~$f_0$ and~$f_{0,1}$. This is simply done by defining:
$$
f_{0,2}(x):=\1e_{B(0,1-\varepsilon)}\tilde p_1(x,T),\quad x \in \Re^3.
$$

\end{enumerate}
\end{Definition}

The last step can be easily justified by the fact that, since~$f_0$ is supported in~$B(0,1-\varepsilon)$:
$$
\norm{f_0-\tilde p_1(\cdot,T)}_{H^1_0(\Re^3)}\geq\norm{f_0-\1e_{B(0,1-\varepsilon)}\tilde p_1(\cdot,T)}_{H^1_0(\Re^3)}.
$$
This scheme is then iterated, constructing a sequence of estimates~$(f_{0,n})_{n\in\Ne^*}$, which may converge to~$f_0$. The following 
theorem, proved in Section~\ref{section_convergence}, ensures that the convergence is geometrical in~$H^1_0(\Re^3)$ under the assumptions:
\begin{itemize}
\item[($i$)] $L$ is the standard wave operator~$\partial_{tt}-\Delta$;
\item[($ii$)] The data are noise-free, {\sl i.e.}~$p_\mathrm{data}=p_\mathrm{exact}$ (see Section~\ref{general_problem}), and the observation
surface~$\S$ is the whole sphere~$S(0,1)$;
\item[($iii$)] The object to reconstruct~$f_0$ is a $C^\infty$ function with compact support in~$B(0,1)$;
\item[($iv$)] $f_{0,1}$ and~$f_{0,2}$ are two consecutive estimates of~$f_0$ obtained as described in Definition~\ref{def_bfn}.
\end{itemize}

\begin{Theorem}\sf
\label{theoreme_principal}
Under the assumptions ($i$)-($iv$), there exists~$1>s>0$ only depending
on~$f_0$,~$\phi$ and the nudging parameter~$k> 0$ such that:
\eqn
\label{inegalite_s}
\norm{f_0-f_{0,2}}_{H^1_0(\Re^3)}\leq s\norm{f_0-f_{0,1}}_{H^1_0(\Re^3)}.
\nqe 
\end{Theorem}

We only consider this ideal framework for the sake of the proof, but numerical results given in Section~\ref{section_illustration} apply with more realistic assumptions and show a better convergence rate.

%%%%%%%%%%%%%%%%%%%%%%%%%%%%%%%%%%%%%%%%%%%%%%%%%%%%%%%%%%%%%%%%%%%%
\section[Proof of Theorem~3]{Proof of Theorem~\ref{theoreme_principal}} %%%%%%%%%%%%%%%%%
\label{section_convergence}   %%%%%%%%%%%%%%%%%%%%%%%%%%%%%%%%%%%%%%%
%%%%%%%%%%%%%%%%%%%%%%%%%%%%%%%%%%%%%%%%%%%%%%%%%%%%%%%%%%%%%%%%%%%%

This section is devoted to the proof of Theorem~\ref{theoreme_principal}. We study the case where~$L=L_0:=\partial_{tt}-\Delta$ is the standard wave operator. Let~$\phi$,~$f_0$ and~$p_{\mathrm{data}}$ be as in Theorem~\ref{theoreme_principal}. In particular, for the sake of the proof, we make the assumptions that the data are noise-free and that, since~$\S=S(0,1)$, the set~$T(\S,\e/2)$ is equal to~$A(1-\e/2,1+\e/2)$. We shall use the notation:
$$
L_k:=\partial_{tt}-\Delta+ k \phi \partial_t,\quad k\in\Re_+.
$$
We first study the existence of the objects used in the BFN approximation scheme.

\subsection{Solution of the backward nudging equation}

The proof of the following theorem can be found in~\cite{bers}, Chapter 2.8.
\begin{Theorem}\sf
Given three~$C^{\infty}$-functions~$\phi$,~$l$ and~$h$ with compact support and~$k\in\Re_+$, 
$T>0$ there exists a unique
$C^{\infty}$ solution~$p$ of the Cauchy problem:
\eqn
\label{eq_nudging_back}
\begin{array}|{ll}.
L_k p(x,t)=-k \phi(x) \partial_t p_{\mathrm{data}}(x,T-t),& (x,t)\in\Re^3\times[0,T], \\
p(x,0)=l(x),& x\in\Re^3, \\
\partial_t p(x,0)=h(x),& x\in\Re^3.
\end{array}
\nqe

Moreover, for every~$(x_0,t_0)\in \Re^3\times\Re_+$,~$p(x_0,t_0)$ only depends on the values 
of~$l$,~$h$,~$\phi$ and~$p_{\mathrm{data}}(\cdot,T-\cdot)$ in the backward characteristic cone:
$$
\norm{x-x_0}\leq\mod{t-t_0}.
$$
\end{Theorem}

In particular, this means that for every sketch~$f_{0,1}$ with compact support and of class 
$C^{\infty}$, the solution~$\tilde p_1$ of the backward nudging equation~\eqref{eq_nudging_back} exists, and for all~$t\in[0,T]$ the 
function~$\tilde p_1(\cdot,t)$ 
has a compact support. This property will be helpful to apply an energy method.

Now, it is easy to show that if~$\tilde p_1$ is solution of~\eqref{eq_nudging_back} with
initial conditions~$l(x)=\tilde p_1(x,0)=p_1(x,T)$ and~$h(x)=\partial_t \tilde p_1(x,0)
=-\partial_t p_1(x,T)$, then the function
$u:=\tilde p_1 - p_{\mathrm{data}}(\cdot,T-\cdot)$ is solution of the Cauchy problem:
\eqn
\label{eq_nudging_reduite}
\begin{array}|{ll}.
L_k u(x,t)=0,&(x,t)\in\Re^3\times[0,T],\\
u(x,0)=p_1(x,T)-p_{\mathrm{data}}(x,T),& x\in\Re^3,\\
\partial_t u(x,0)=-\partial_t(p_1-p_{\mathrm{data}})(x,T),& x\in\Re^3,
\end{array}
\nqe
and satisfies:
$$
E_u(0)=E_{p_1-p_{\mathrm{data}}}(T)=\dfrac{1}{2}\norm{f_{0,1}-f_0}_{H^1(\Re^3)}^2.
$$
Since, from the definition of~$f_{0,2}$:
$$
E_u(T)\geq\dfrac{1}{2}\norm{\tilde p_1(\cdot,T)-f_0}_{H^1(\Re^3)}^2\geq
\dfrac{1}{2}\norm{f_{0,2}-f_0}_{H^1(\Re^3)}^2,
$$
we only need to prove that:
$$
E_u(T)\leq s E_u(0),\quad \hbox{for some }0<s<1,
$$
to complete the proof of Theorem~\ref{theoreme_principal}. This is the purpose of the next section.

%%%%%%%%%%%%%%%%%%%%%%%%%%%%%%%%%%%%%%%%%%%%%%%%%%%%%%%%%%%%%%%%%%%%%%%%%
\subsection{Energy inequality}

We have the following identity for every solution~$u$ of~\eqref{eq_nudging_back}:
$$
\partial_t u L_k u=0=\partial_t \left [ \dfrac{\norm{\nabla u}^2}{2}+\dfrac{(\partial_t u)^2}{2} \right ]
+k \phi (\partial_t u)^2-\sum_{i=1}^{3}\partial_{x_i}(\partial_t u\, \partial_{x_i}u).
$$
Recall that~$u$ has a compact support for every~$t$, thus, integrating the latter
equality over the whole space and for~$0<t<T$, we get:
\eqn
\label{egalite_energie}
E_u(T)-E_u(0)=-k\int_0^T\int_{\Re^3}\phi(\partial_t u)^2.
\nqe
Consequently, the loss of energy of~$u$ is exactly the amount of its 'energy' passing through the support of~$\phi$ during time~$T$.

In \cite{finch}, in order to invert the spherical Radon transform, Finch \textit{et al}
proved the following \emph{trace identity}, which will be usefull to estimate \eqref{egalite_energie}:
 
\begin{Theorem}\sf
Suppose~$h_i\in C_0^{\infty}(B(0,\rho))$ and~$u_i$ is the solution of \eqref{eq_onde_exemple}
for~$l=0$ and~$h=h_i$,~$i=1,2$. Then we have the identity
$$
\dfrac{1}{2}\int_{\Re^3}h_1h_2=\dfrac{-1}{\rho}\int_0^{\infty}\int_{S(0,\rho)}tu_1\partial_{tt}u_{2}.
$$
\end{Theorem}
 
From this we will obtain a key estimate of the 'energy' passing through the sphere~$S(0,\rho)$. Indeed, if~$u$ is solution of \eqref{eq_onde_exemple} with~$h=0$, taking~$h_1=f$ and~$h_2=\Delta f$ in the theorem yields:
$$
\dfrac{1}{2}\int_{\Re^3}f\Delta f=\dfrac{-1}{\rho}\int_0^{\infty}\int_{S(0,\rho)}t(\partial_t u)^2,
$$
which gives:
\eqn
\label{estimation_sphere}
\dfrac{1}{2}\norm{f}_{H^1_0(\Re^3)}^2=\dfrac{1}{\rho}\int_0^{\infty}\int_{S(0,\rho)}t(\partial_t u)^2.
\nqe

Let us introduce the solution~$u_0$ of the Cauchy problem:
$$
\begin{array}|{ll}.
L_0 u_0(x,t)=0,&(x,t)\in\Re^3\times[0,T],\\
u_0(x,0)=p_1(x,T)-p_{\mathrm{data}}(x,T),& x\in\Re^3,\\
u_{0t}(x,0)=-\partial_t p_1 (x,T)+\partial_t p_{\mathrm{data}}(x,0),& x\in\Re^3,
\end{array}
$$
which is nothing but~$(p_1-p_{\mathrm{data}})(\cdot,T-\cdot)$.

Then~$u=u_0+v$ where~$v$ is solution of:
$$
\begin{array}|{ll}.
L_0 v(x,t)=-k\phi\partial_t(u_0+v)(x,t),&(x,t)\in\Re^3\times[0,T],\\
v(x,0)=\partial_t v(x,0)=0,&x\in\Re^3.
\end{array}
$$

The calculation of the energy of~$v$ yields:
$$
E_v(T)=-k\int_0^T\int_{\Re^3}\phi \partial_t u_{0}\partial_t v-k\int_0^T\int_{\Re^3}\phi (\partial_t v)^2.
$$
Thanks to the Cauchy-Schwarz inequality, one has:
\eqn
\label{en_v}
E_v(T)\leq k\norm{\sqrt{\phi}\partial_t u_{0}}_{L^2([0,T]\times
\Re^3)}H(T)-kH(T)^2,
\nqe
where:
$$
H(T):=\norm{\sqrt{\phi}\partial_t v}_{L^2([0,T]\times\Re^3)}.
$$
Moreover, we have:
$$
E_v(T)\geq\int_{\Re^3}(\partial_t v)^2(\cdot,T)\geq\int_{\Re^3}\phi (\partial_t v)^2(\cdot,T)=\dfrac{\ud}{\ud T}\left (H^2\right )(T).
$$
So that, dividing by~$H(T)$ in \eqref{en_v} we find:
$$
\dfrac{\ud}{\ud T}H(T)+kH(T)\leq k\norm{\sqrt{\phi}\partial_t  u_0}_{L^2([0,T]\times\Re^3)}.
$$
A classical calculation shows that this differential inequation yields:
$$
H(T)\leq \int_0^T k\norm{\sqrt{\phi}\partial_t u_{0}}_{L^2([0,s]\times\Re^3)}
e^{k(s-T)}\ud s,
$$
and since the map~$s\mapsto\norm{\sqrt{\phi}\partial_t u_{0}}_{L^2([0,s]\times\Re^3)}$ is increasing with~$s$, we finally get:
$$
H(T)\leq \norm{\sqrt{\phi}\partial_t u_{0}}_{L^2([0,T]\times\Re^3)}
(1-e^{-kT}).
$$
This inequality expresses the fact that despite the feedback part of Equation~\eqref{eq_nudging_back}, its solution~$u$ keeps a certain 
amount of energy on the
support of~$\phi$, which is proportional to the same energy for~$u_0$. Indeed:
\eqa
\nonumber \norm{\sqrt{\phi}\partial_t u}_{L^2([0,T]\times\Re^3)}
&\geq &\norm{\sqrt{\phi}\partial_t u_{0}}_{L^2([0,T]\times\Re^3)}
-H(T)\\
\nonumber &\geq& \norm{\sqrt{\phi}\partial_t u_{0}}_{L^2([0,T]\times\Re^3)}e^{-kT},
\aqe
thus:
\eqn
\label{estimation_u}
\int_0^T\int_{\Re^3}\phi (\partial_t u)^2\geq \underset{A}{\underbrace{\left ( \int_0^T\int_{\Re^3}\phi (\partial_t u_{0})^2 \right )}}e^{-2kT}.
\nqe
The last step of this proof is the estimation of~$A$. We have:
$$
A\geq \int_{A(1-\varepsilon/2,1+\varepsilon/2)}\int_0^T 
(\partial_t u_{0})^2=\int_{A(1-\varepsilon/2,1+\varepsilon/2)}\int_0^T (\partial_t p_1-\partial_t p_{\mathrm{data}})^2(T-t),
$$
so that:
$$
A\geq \int_{A(1-\varepsilon/2,1+\varepsilon/2)}\int_0^{T} (\partial_t p_1-\partial_t p_{\mathrm{data}})^2
=\int_{1-\varepsilon/2}^{1+\varepsilon/2}\int_0^{T}\int_{S(0,\rho)}
(\partial_t p_1-\partial_t p_{\mathrm{data}})^2.
$$
Since~$p_1 -p_{\mathrm{data}}$ vanishes on~$A(1-\varepsilon/2,1+\varepsilon/2)$ for every~$t>T$:
$$
A\geq \dfrac{1}{T}\int_{1-\varepsilon/2}^{1+\varepsilon/2}\int_0^{\infty}\int_{S(0,\rho)}
t(\partial_t p_1-\partial_t p_{\mathrm{data}})^2.
$$
Finally, applying~\eqref{estimation_sphere} to~$p_1-p_{\mathrm{data}}$, one has:
$$
A\geq \dfrac{1}{T}\int_{1-\varepsilon/2}^{1+\varepsilon/2}\rho \dfrac{1}{2}\norm{f-f_{0,1}}_{H^1_0(\Re^3)}^2 \ud \rho,
$$
which yields:
\eqn
\label{estimate_A}
A\geq\dfrac{\varepsilon}{T}E_u(0).
\nqe
Combining the above estimate with~\eqref{egalite_energie} and~\eqref{estimation_u}, we get:
$$
E_u(T)\leq (1-\dfrac{k\varepsilon}{T}e^{-2kT})E_u(0),
$$
which was the inequality to prove, with~$s:=1-\frac{k\varepsilon}{T}e^{-2kT}<1$.

\begin{Remark}
Since we know that~$T=2$, there exists an optimal choice for the nudging parameter~$k$
which leads to~$s=1-\frac{\varepsilon}{8e}$. Even though this theoretical convergence rate
does not seem satisfactory, we shall see in Section~\ref{section_illustration} that the
actual convergence is much faster in practice.
\end{Remark}

%%%%%%%%%%%%%%%%%%%%%%%%%%%%%%%%%%%%%%
\section{Variational approach}%%%%%%%%
%%%%%%%%%%%%%%%%%%%%%%%%%%%%%%%%%%%%%%
\label{var}

We deal now with a variational formulation of the inverse TAT problem, that leads to some new reconstruction techniques.

In this section, we consider the reconstruction of an object~$f_0\in H^1_0(B(0,1))$ from a set of observations~$p_\mathrm{data}=p_\mathrm{exact}+p_\mathrm{noise}$. Here,~$p_\mathrm{exact}:=\W f_0$ where the linear operator~$\W$ maps an initial condition~$f$ to its related solution of the Cauchy problem~\eqref{eq_onde_exemple} with a null derivative initial condition~$h=0$, restricted to a non-empty, open observation set~$A_\e\subset A(1-\e,1+\e)$. Namely,~$\W$ can be written~$\W=\Phi W$, with:
$$
\operator{W}{H_0^1(B(0,1))}{C^0([0,T];H^1_0(B(0,R)))}{f}{Wf,}
$$
where~$Wf$ is the weak solution of~\eqref{eq_onde_exemple} with the initial conditions~$l=f$ and~$h=0$.
Even though~$Wf$ is defined on the whole space~$\Re^3$, we can choose~$R>0$ large enough for any solution of~\eqref{eq_onde_exemple} to keep a null trace on~$S(0,R)$, such that the restriction of~$W f$ to~$B(0,R)$ is in~$H^1_0(B(0,R))$ 
for every~$t\in[0,T]$. The operator~$\Phi$ is the restriction operator:
$$
\operator{\Phi}{L^2([0,T]\times\Re^3)}{L^2([0,T]\times A_\e)}{g}{g_{|_{A_\e}}.}
$$
Our purpose is to solve the inverse problem: $$\W f=p_\mathrm{data},$$
by means of the following minimization problem:
$$
(\P_\alpha) \left|
\begin{array}{rl}
\hbox{Minimize}& J_\alpha(f)\\
\hbox{s.t.}&f\in H^1_0(B(0,1)).
\end{array}
\right.
$$
Here,~$J_\alpha$ is the Tikhonov functional:
$$
\operator{J_\alpha}{H^1_0(B(0,1))}{\Re_+}{f}{\frac{1}{2}\int_0^T \norm{\W f-p_\mathrm{data}}^2_{L^2(A_\varepsilon)}\ud t+\frac{\alpha}{2}\norm{f}^2_{L^2(B(0,1))},}
$$
with~$\alpha>0$. Since~$\W$ is a bounded linear operator,~$J_\alpha$ is a strictly convex (thanks to the regularization parameter~$\alpha$), differentiable and coercive functional, and any solution of~$(\P_\alpha)$ is characterized by:
$$
\nabla J_\alpha(f)=\W^*(\W f-p_\mathrm{data})+\alpha f=0.
$$
For any~$\alpha>0$, the above equation has a unique solution:
$$
f_\alpha=(\W^*\W+\alpha I)^{-1}\W^* p_\mathrm{data},
$$
so that Problem~$(\P_\alpha)$ is well-posed. We have chosen to solve this latter by means of the conjugate gradient method, which requires, at each iteration, the computation of~$\nabla J_\alpha$. 

Now, one has, for every~$\psi \in C^\infty_0(B(0,1))$ and~$f\in H_1^0(B(0,1))$:
\begin{eqnarray*}
\scal{\nabla J_\alpha(f)}{\psi}_{L^2(B(0,1))}&=&\scal{\W^*(\W f-p_\mathrm{data})}{\psi}+\alpha\scal{f}{\psi}\\
&=&\underbrace{\scal{\W f-p_\mathrm{data}}{\W \psi}}_{G_{f,\psi}}+\alpha\scal{f}{\psi},
\end{eqnarray*}
with:
\begin{eqnarray*}
G_{f,\psi}&=&\int_0^T\int_{A_\e}(Wf-p_\mathrm{data})W\psi\\
&=&\int_0^T\int_{B(0,R)}\1e_{A_\e}(Wf-p_\mathrm{data})W\psi.\\
\end{eqnarray*}
Considering the weak solution~$u^*$ of the adjoint equation:
\eqn
\label{eq_adj}
\begin{array}|{l}.
\partial_{tt}u^*(x,t)-\Delta u^*(x,t)=\1e_{A_\e}(Wf-p_\mathrm{data}),\\
u^*(x,T)=0,\\
\partial_{t}u^*(x,T)=0,
\end{array}
\nqe
and since~$W\psi\in C^\infty_0(B(0,R))$, we have:
$$
G_{f,\psi}=\int_0^T\int_{B(0,R)}\partial_{tt} u^* W\psi+\int_0^T\int_{B(0,R)}\scal{\nabla u^*}{\nabla W\psi}.
$$
Since~$B(0,R)$ has been chosen large enough, applying Green's formula and integrating by parts yield:
\begin{eqnarray*}
G_{f,\psi}&=&\int_0^T\int_{B(0,R)}u^* L W\psi-\int_{B(0,R)}\partial_t u^* (x,0)W\psi(0)\ud x\\
&=&\scal{-\partial_t u^*(\cdot,0)}{\psi}.
\end{eqnarray*}
Finally:
$$
\nabla J_\alpha(f)=-\partial_t u^*(\cdot,0)+\alpha f.
$$

This calculation ensures that each iteration of the conjugate gradient algorithm necessitates successively the 
computations of~$W f$ and of its adjoint state~$u^*$, whose algorithmic complexity is comparable to one iteration
of the BFN algorithm. This should be noted that the computation of the adjoint state requires the storage of~$Wf - p_{\mathrm{data}}$ on~$A_\e$, which is not required in the BFN algorithm.

Note that even if we have formulated the calculation of~$\nabla J_\alpha$ in a continuous framework, its
discrete counterpart leads to similar considerations. Obviously, the behaviour of the method strongly depends
on the scheme chosen to compute solutions of~\eqref{eq_onde_exemple}.

%%%%%%%%%%%%%%%%%%%%%%%%%%%%%%%%%%%%%%%%%%%%%%%%%%%%%%%%%%%%%%%%%%%%%%%%%%%
\section{Numerical implementation}%%%%%%%%%%%%%%%%%%%%%%%%%%%%%%%%%%%%%%
%%%%%%%%%%%%%%%%%%%%%%%%%%%%%%%%%%%%%%%%%%%%%%%%%%%%%%%%%%%%%%%%%%%%%%%%%%%%
\label{section_illustration}
In order to test the numerical behavior of the BFN method, we have implemented it in several situations.
In particular, we have considered the case where the data are available only on a half-sphere, so that the
illustration fits the practical experimental set-up.

%%%%%%%%%%%%%%%%%%%%%%%%%%%%%%%%%%%%%%%%%%%%%%%%%%%%%%%%%%%%%
\subsection{Implementation framework}%%%%%%%%%%%%%%%%%%%%%%%%
%%%%%%%%%%%%%%%%%%%%%%%%%%%%%%%%%%%%%%%%%%%%%%%%%%%%%%%%%%%%%
\label{implementation}

We now describe the numerical context chosen to compare Back and Forth Nudging (BFN) and the Conjugate Gradient (CG) techniques to a reasonable algorithm, namely the Neumann Series (NS) (since it is shown in~\cite{stefanov_algorithm} that the use of NS significantly improves the Time Reversal (TR) algorithm, we mainly chose to consider NS results).

We work now in a 2D framework for numerical costs reasons. The object to reconstruct is in~$[-0.5,0.5]^2$ gridded with~$256^2$ pixels, such that the space step is~$\delta x=1/256$.

Sensors are located on the circle~$S(0,\sqrt2/2)$ surrounding the object, on which various distributions are considered. If no additional information is given, 800 detectors are used when the data are complete. Nevertheless some reconstructions were computed with less detectors in order to bring to light the robustness of the techniques. When a multiplicative white Gaussian noise is added to the data, their level is set to~$15\%$.

The final time~$T=\sqrt2$ is then taken minimal in order to leave waves traveling to the sensors when sound speed is constant with value~$1$. Note that the four studied methods can significantly be improved by increasing~$T$ (\textit{cf.}~\cite{stefanov_algorithm}). We use the domain~$[-0.5-\sqrt2-\e;0.5+\sqrt2+\e]^2$ for the computation, where~$\e$ represents a few pixels, which is large enough to avoid reflection effects on its boundary.

The model is computed by means of the classical finite difference time domain (FDTD) method (explicit Euler scheme and classical five points discretization of the Laplace operator~$\Delta_{\delta x}$).

\begin{figure}[!ht]
\centering
\includegraphics[width=.3\textwidth,height=.3\textwidth]{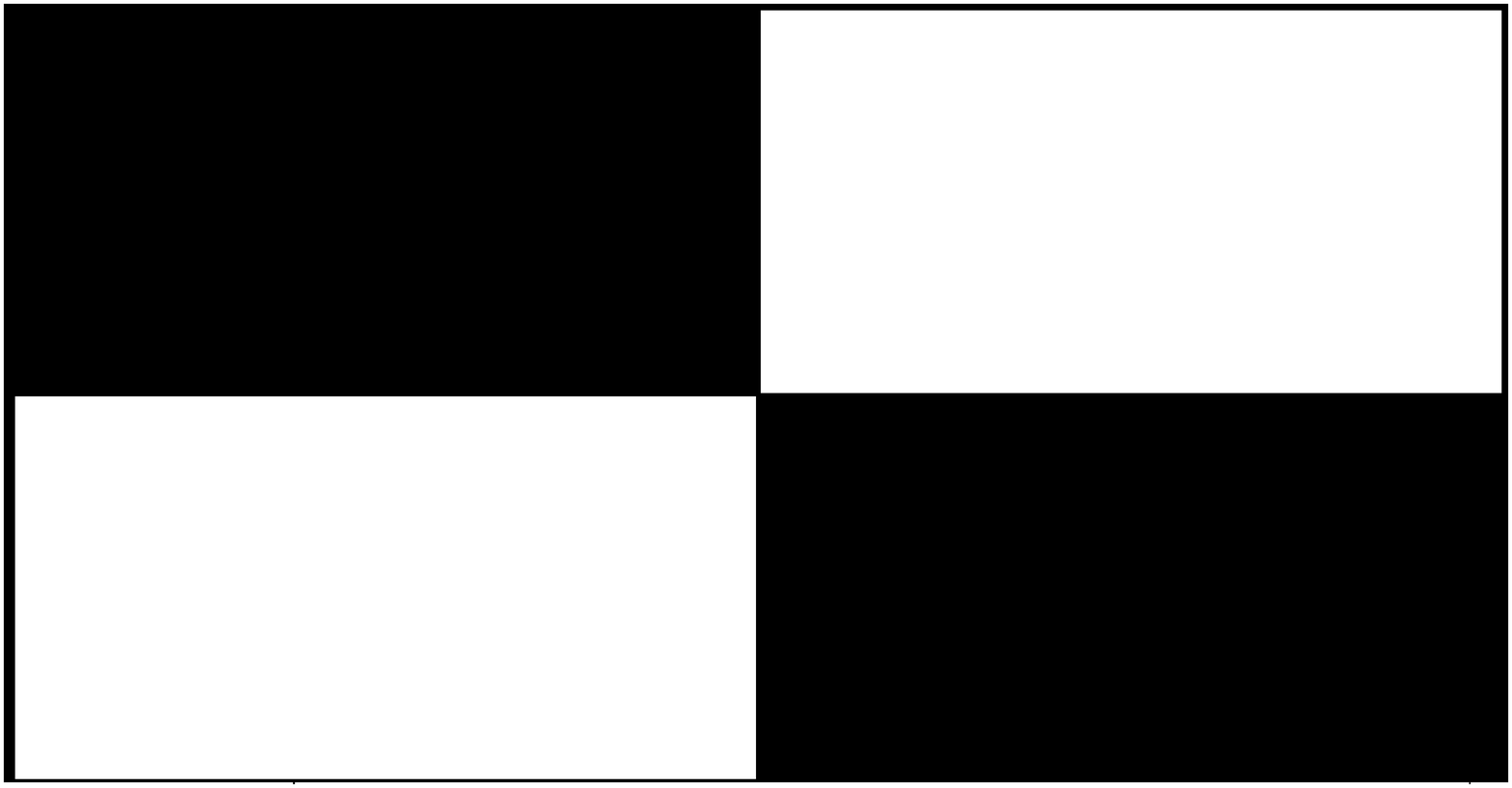}
\includegraphics[width=.3\textwidth,height=.3\textwidth]{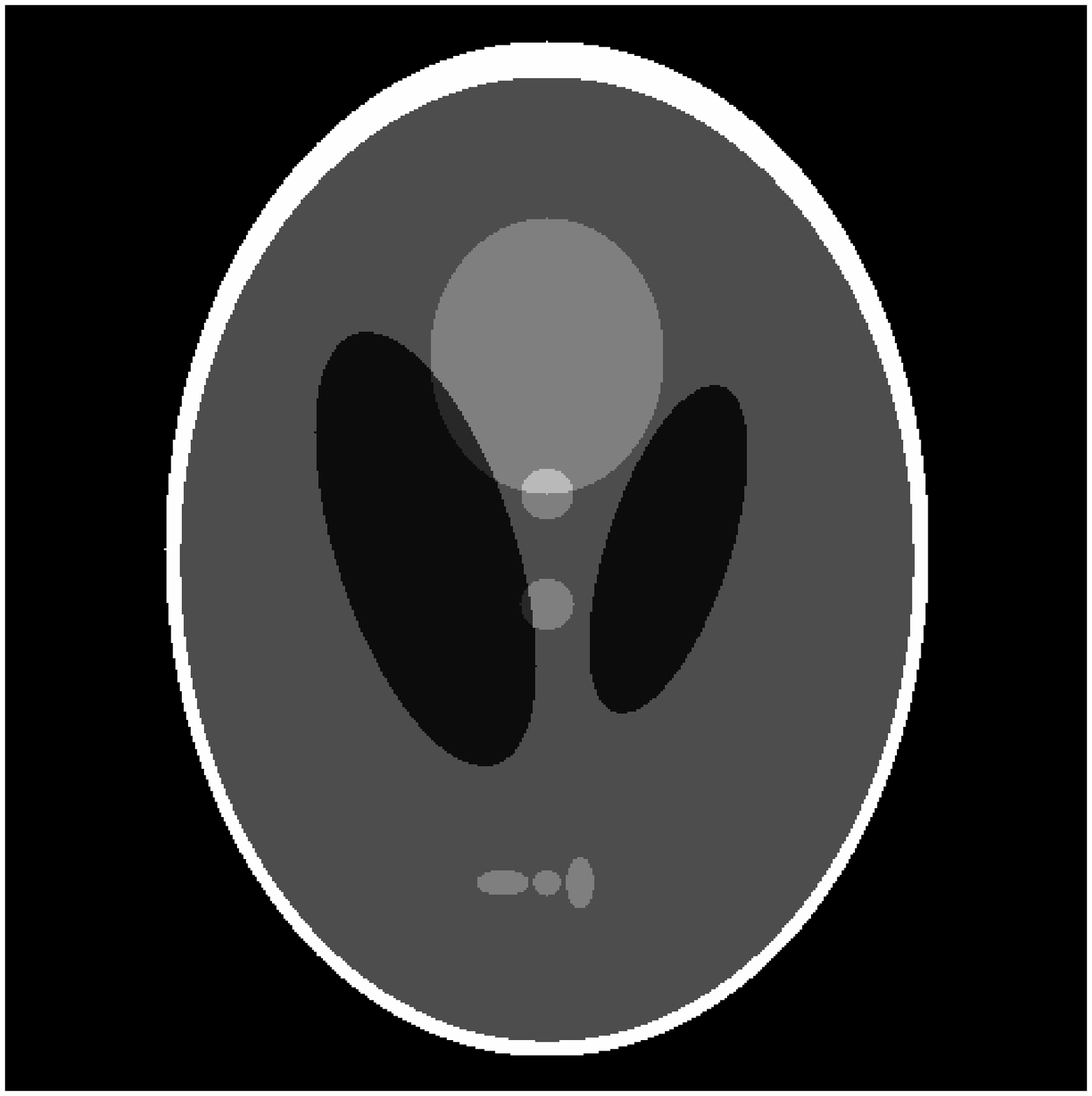}
\includegraphics[width=.3\textwidth,height=.3\textwidth]{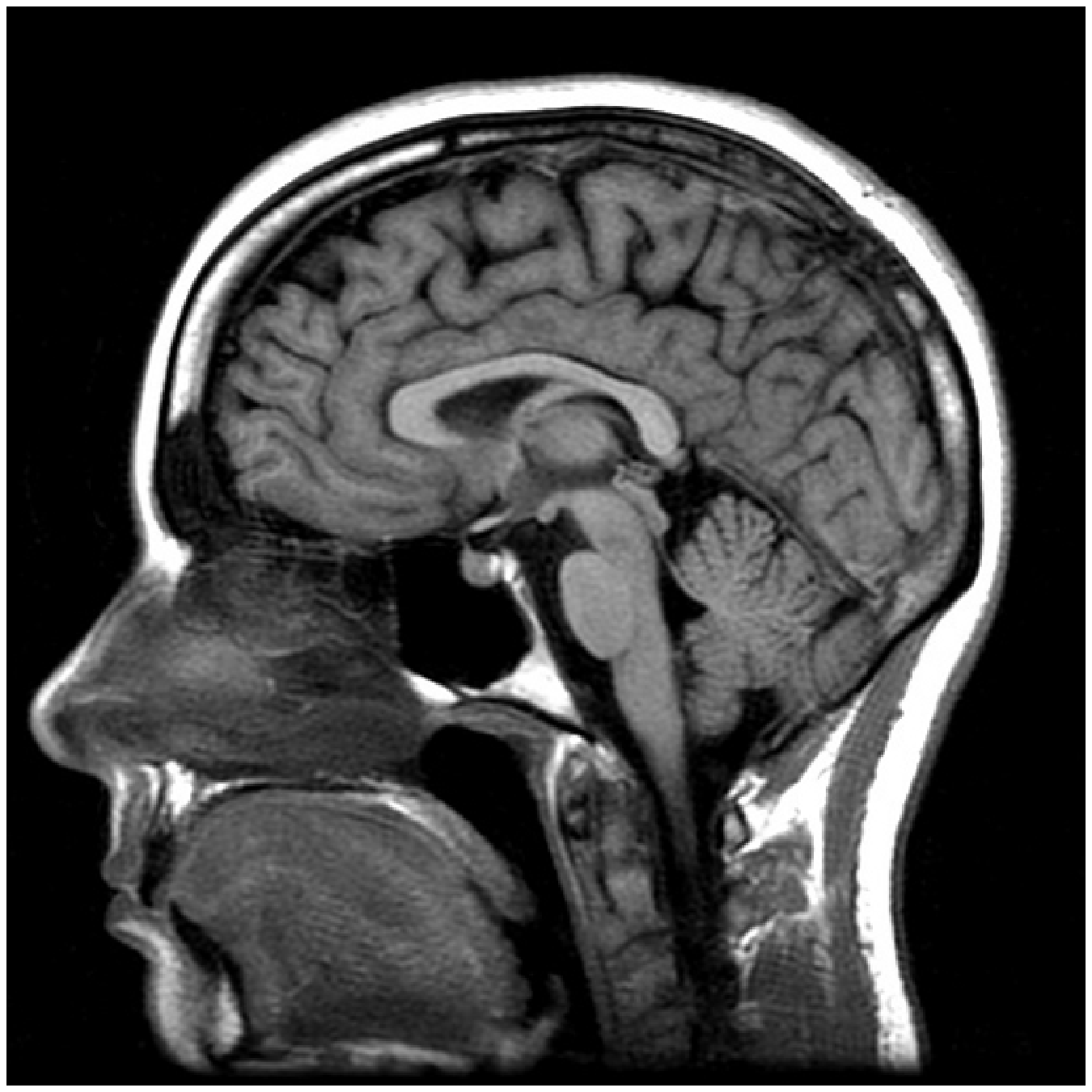}
\caption{\label{objets}The initial objects $f_0$: squares, Shepp-Logan and skull.}
\end{figure}

Figure~\ref{objets} represents the 2-dimensional objects to reconstruct, latter denoted by~$f_0$. Squares able to give rise to typical phenomenons of observation of wave-like systems while the second object, known as the Shepp-Logan phantom, allows good analysis and comparisons of the ability of the methods. The third object has been considered for combining sharp contours with smoother areas. Each is implemented as a 256 by 256 matrix.

\begin{figure}[!ht]
\centering
\includegraphics[width=.3\textwidth,height=.3\textwidth]{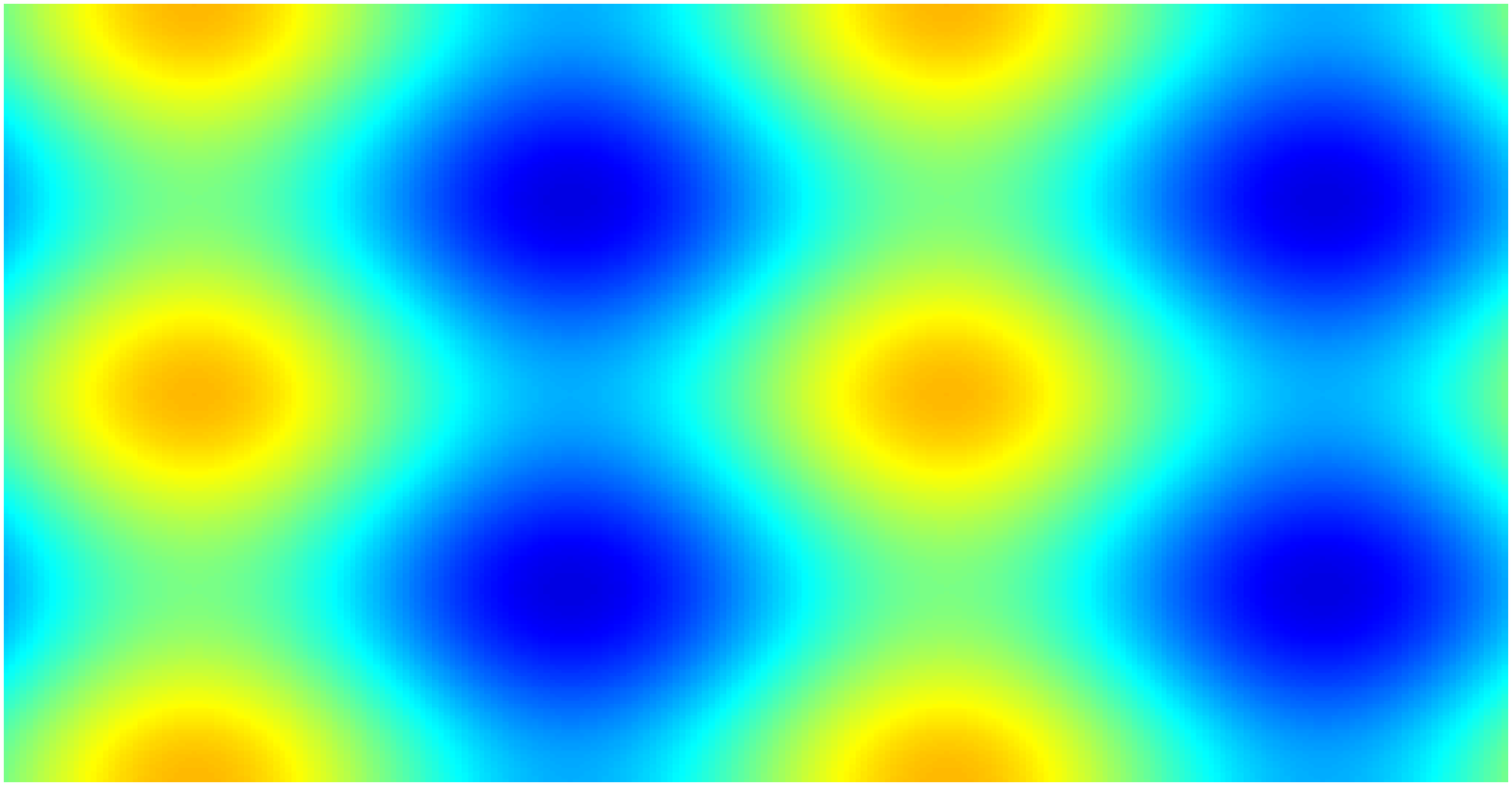}
\includegraphics[width=.3\textwidth,height=.3\textwidth]{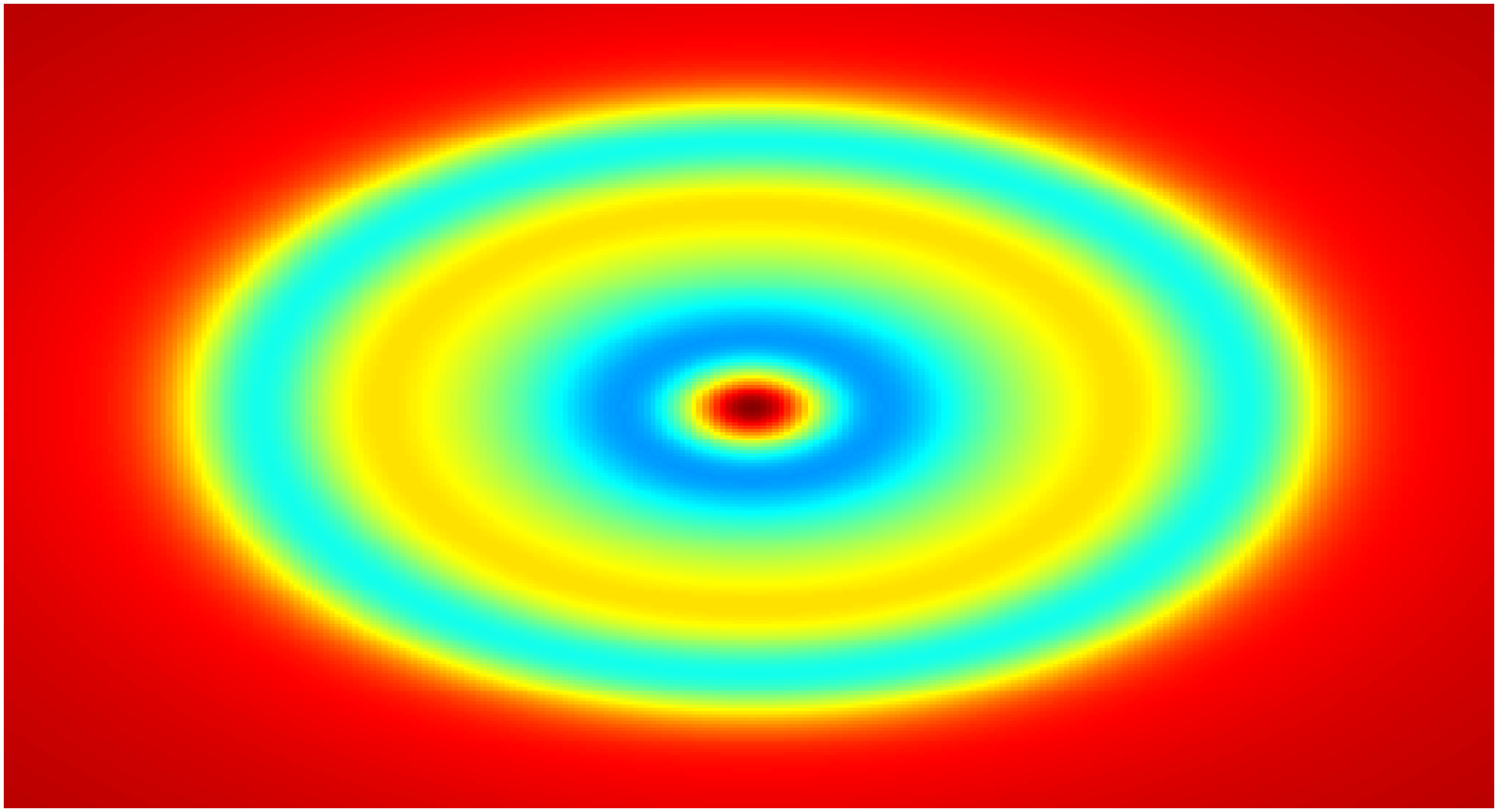}
\includegraphics[width=.3\textwidth,height=.3\textwidth]{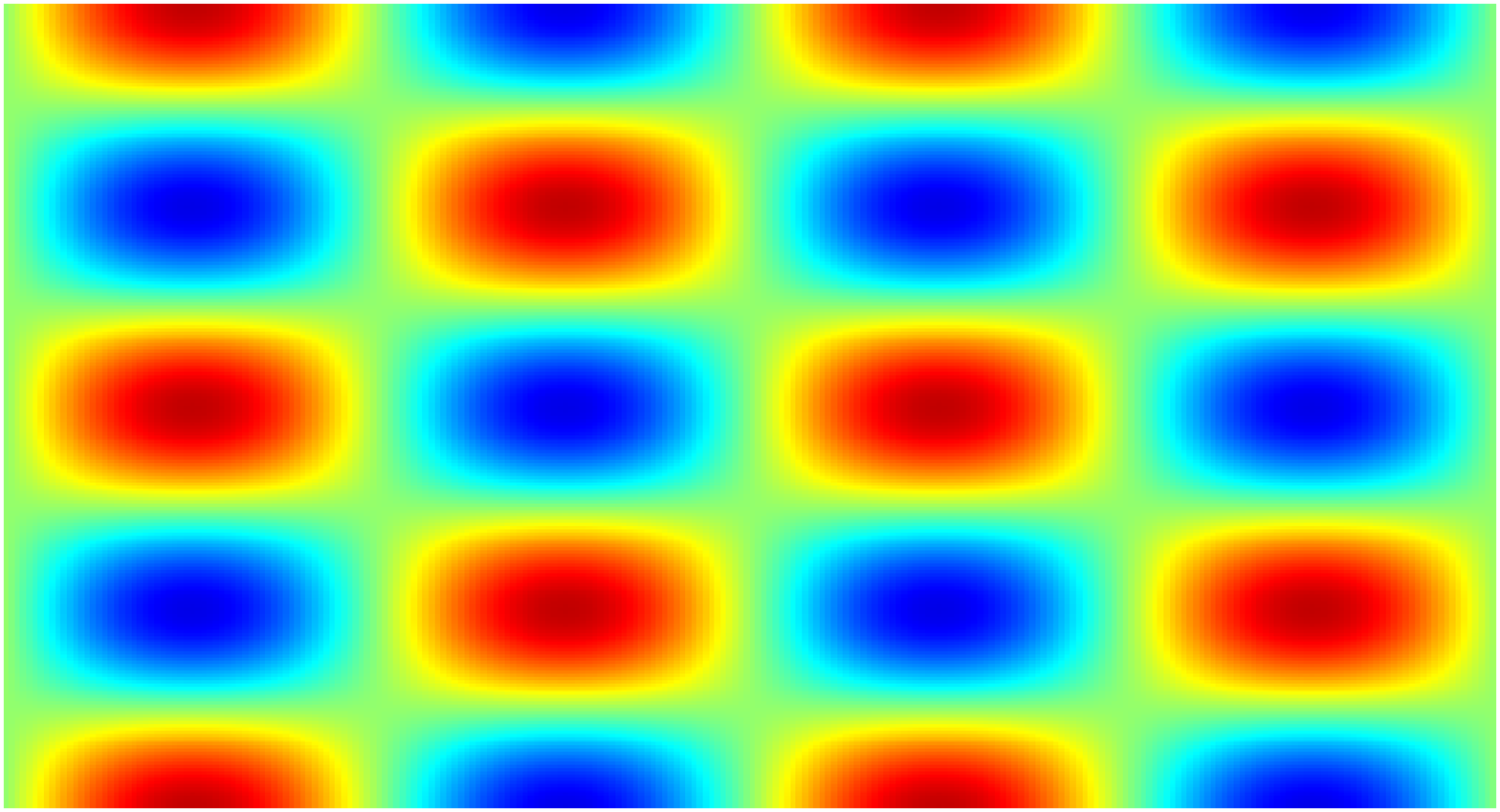}
\caption{\label{speeds}The different variable speeds: NTS, TS1 and TS2 (from equations~\eqref{eqvit}).}
\end{figure}

Four sound speeds are used. First a constant one~$c=1$, then three other one, variable, which are chosen as in~\cite{stefanov_algorithm}, respectively defined as:
\begin{eqnarray}
\label{eqvit}
c(x,y)&=&1+0.2\sin(2\pi x)+0.1\cos(2\pi y),\nonumber\\
c(x,y)&=&\frac{9(x^2+y^2)}{1+9(x^2+y^2)}+\exp(-90(x^2+y^2))\\
&&-0.4\exp(-10(3\sqrt{x^2+y^2}-2)^2),\nonumber\\
c(x,y)&=&1.25+\sin(2\pi x)\cos(2\pi y).\nonumber
\end{eqnarray}
They are plotted on figure~\ref{speeds} and abbreviated respectively NTS, TS1 and TS2, as the first is a Non-Trapping Speed while the two other are Trapping Speeds. No cutoff is used to smooth the transition between variable internal speed and constant external speed (always being~1).

Concerning the parameters of the algorithms, the regularization parameter~$\alpha$ in CG method equals~0. Other values were tested that stabilized the algorithm but always blur the reconstructions so that we get higher errors than without regularization. The observed convergence rates of CG and NS methods are fairly high, so that only~$10$ iterations yield a good approximation of the best achievable solutions. With regard to the two nudging coefficients, while we can theoretically affirm that the larger it is, the faster the solution tends to the observations (thanks to the energy equality~\eqref{egalite_energie}), the discretization of the time-derivative of the solution imposes a coefficient~$k$ not larger than~$1/\delta t$. In order to obtain the best convergence and to make the choice of~$k$ less arbitrary (and fitting with the limits of the numerical method), it is defined as $0.9/\delta t$.

The following tables give the relative mean square error obtained with the best estimation obtained after at most~10 iterations and, in brackets, the corresponding number of iterations required. Note that BFN, CG and NS furnishes convex error plots with a minimum.

Finally, as it has been suggested in~\cite{bdch}, we also have used an artificial numerical attenuation to offset spurious high-frequency effects and noise, with interesting consequences. Then for BFN, NS and CG techniques, the numerical model becomes:
$$
\frac{u_{n+1}-2u_n+u_{n-1}}{\delta t^2}=\Delta_{\delta x}u_n+\delta x^{\e}\Delta_{\delta x}\frac{u_n-u_{n-1}}{\delta t},
$$
where the value for~$\e$ is empirically set and an optimal value appears for each method, which are approximately~2 for BFN and~1.7 for CG and NS in a complete data situation.

The attenuation being only introduced in a numerical improvement view, it is set to attenuate both back and forth solutions, should the case arise, and when the data 
are corrupted by noise.

%%%%%%%%%%%%%%%%%%%%%%%%%%%%%%%%%%%%%%%%%%%%%%%%%%%%%%%%%%%%%%%%%%
\subsection{Implementation results}%%%%%%%%%%%%%%%%%%%%%%%%%%%%%%%
%%%%%%%%%%%%%%%%%%%%%%%%%%%%%%%%%%%%%%%%%%%%%%%%%%%%%%%%%%%%%%%%%%
\label{numres}

The first results are given in Table~\ref{first} in a noiseless complete data situation. One sees that variable speeds can have various effects, depending on the method we use. NS offers the best reconstruction, then BFN and finally CG.

\begin{table}[!h]\centering
$$
\begin{array}.{lccc}.
\text{Speed map}&\text{BFN}&\text{CG}&\text{NS}\\
c\equiv1&2.3 (10)&7.1 (3)&1.4 (10)\\
\text{NTS}&3.2 (10)&7.5 (4)&2.3 (10)\\
\text{TS1}&7.1 (10)&16.5 (5)&1.6 (10)\\
\text{TS2}&7 (10)&10.7 (5)&6.6 (10)
\end{array}
$$
\caption{\label{first}Noiseless complete data for Shepp-Logan phantom: relative mean square error (number of iterations, limited to~10).}
\end{table}

These techniques can be consequently improved by increasing the number of iterations as we can see in Table~\ref{cent} for the BFN, where~100 iterations instead of~10 are computed.

\begin{table}[!h]\centering
$$
\begin{array}.{cccc}.
c\equiv1&\text{NTS}&\text{TS1}&\text{TS2}\\
1&1.5 (100)&1.2 (100)&4.1 (100)
\end{array}
$$
\caption{\label{cent}Noiseless complete data from Shepp-Logan phantom: relative mean square error (BFN,~100 iterations).}
\end{table}

%Noisy complete data results on Table~\ref{noise}. One notices~:~$15\%$ level noise has globally more effects than perturbations of constant speed (cf \textit{e.g.} the different errors for BFN depending on the speeds, with or without noise), TR good first rough estimate (so as NS iterations can be useless), BFN and CG more sensible to noise than to variable speeds in comparison with TR, which reacts worse to variable speed and keeps robust to noise).
%
%\begin{table}[!h]\centering
%$$
%\begin{array}.{lcccc}.
%\text{Speed map}&\text{BFN}&\text{TR}&\text{CG}&\text{NS}\\
%c\equiv1&28.8 (3)&26.8&35.5 (8)&26.8 (1)\\
%%\text{NTS}&27 (3)&24.8&33.3 (9)&24.8 (1)\\
%%\text{TS1}&38.8 (2)&42.8&37.7 (6)&35.8 (2)\\
%%\text{TS2}&28.9 (3)&27.4&34.3 (5)&27.4 (1)
%\end{array}
%$$ 
%\caption{\label{noise}Noisy complete data ($15\%$ noise).}
%\end{table}
%
%Then same situation, with constant speed, but trying to avoid noise effects with an artificial numerical attenuation. 

%Tests to optimize~$\e$ for each method on Table~\ref{besteps}. Obvious minimum. Depends highly on settings (position, repartition and number of sensors).
%
%\begin{table}[!h]\centering
%$$
%\begin{array}.{lccccc}.
%&\text{BFN}&\text{TR}&\text{CG}&\text{NS}\\
%\e=1.6&25.5 (50)&20.4&21 (1)&17.6 (3)\\
%\e=1.7&22.9 (48)&18.4&19.8 (1)&17.1 (2)\\
%\e=1.8&20.6 (53)&17.3&20.4 (1)&17.3 (1)\\
%\e=1.9&18.8 (53)&17.3&23.7 (1)&17.3 (1)\\
%\e=2&17.6(48)&18.3&29.9 (1)&18.3 (1)\\
%\e=2.1&17.4(52)&&34.9 (2)&\\
%\e=2.2&18.1 (52)&&35.9 (2)&
%\end{array}
%$$ 
%\caption{\label{besteps}Constant speed,~$15\%$ noise, complete data.}
%\end{table}

When adding a~$15\%$ level noise in the constant speed situation, one sees on Table~\ref{att} that any method offers an interesting reconstruction: one gets less than~$20\%$ error for~$15\%$ noise thanks to the attenuating scheme.

\begin{table}[!h]\centering
$$
\begin{array}.{ccc}.
\text{BFN}&\text{CG}&\text{NS}\\
%\e=0&28.8 (3)&35.5 (8)&26.8 (1)\\
17.6 (10)&19.8 (2)&17.1 (2)
\end{array}
$$ 
\caption{\label{att}Noisy complete data ($15\%$ noise) from Shepp-Logan phantom with the constant speed. Reconstruction with attenuating sheme: relative mean square error (number of iterations, limited to 10).}
\end{table}

\begin{figure}[!h]
\centering
\includegraphics[width=.4\textwidth,height=.4\textwidth]{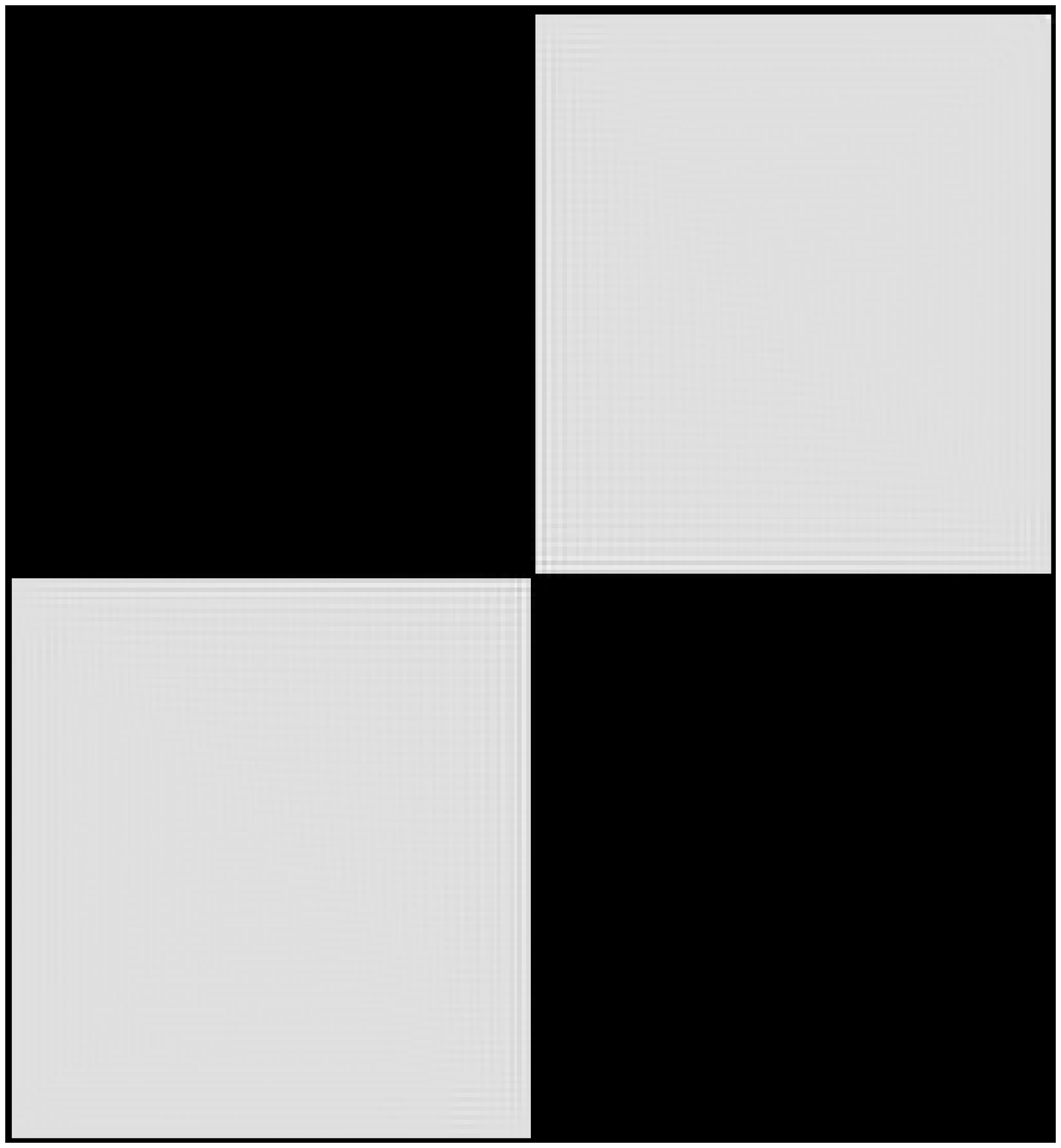}
\includegraphics[width=.395\textwidth,height=.395\textwidth]{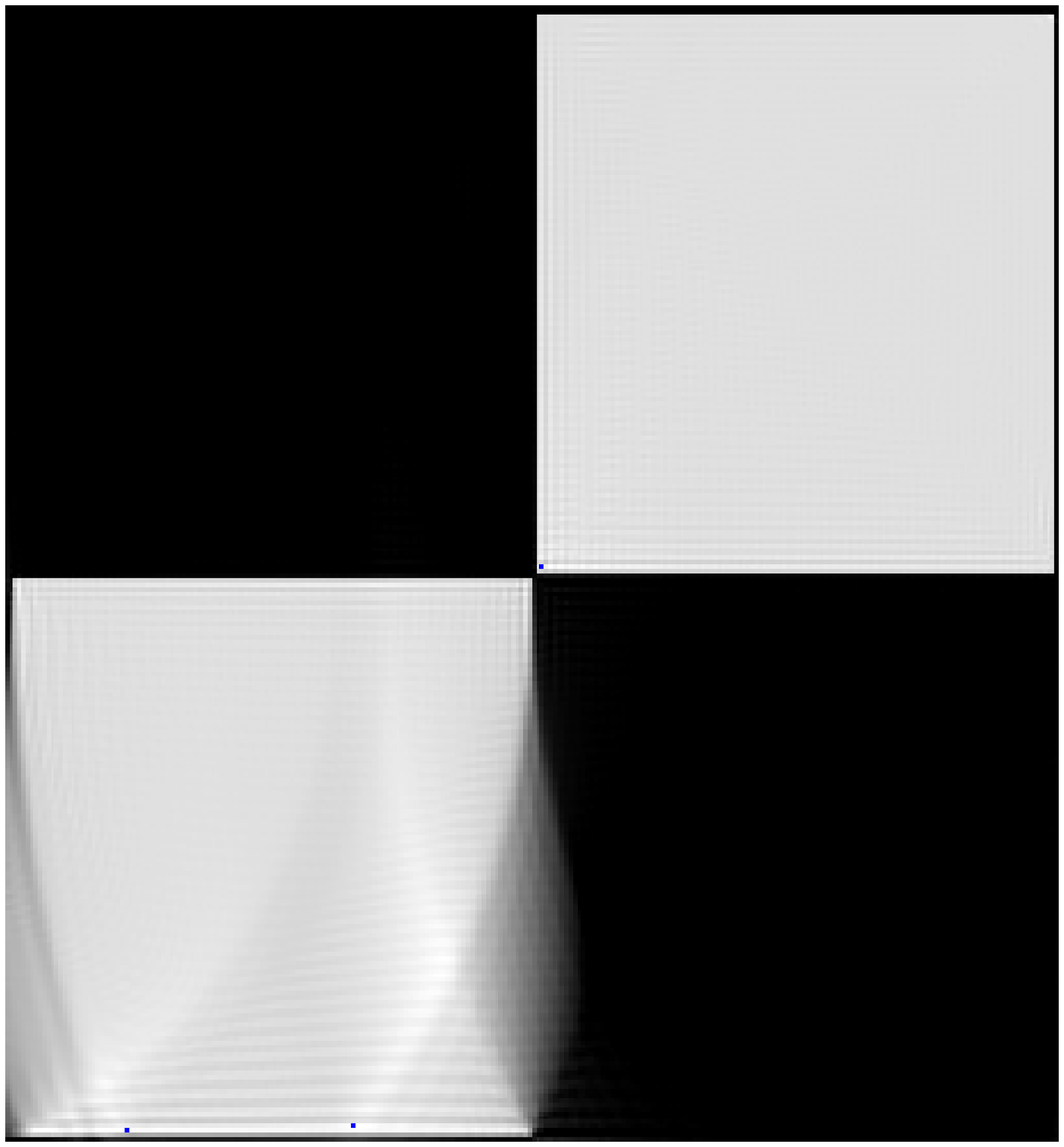}
\caption{\label{incarre}Reconstructions by~100 BFN iterations. Squares observed singularities. Complete ($1.1\%$ error) and half-circle ($10.4\%$) noiseless data with constant speed. }
\end{figure}

About the propagation of the singularities and their observation, we put the light on the effect of getting incomplete data (in a noiseless situation) in Fig.~\ref{incarre}. The sensors are located on the upper half circle, that illustrate the practical case where the observed body cannot be surrounded by the observation surface.

\begin{figure}[!h]\centering
\includegraphics[width=.5\textwidth,height=.5\textwidth]{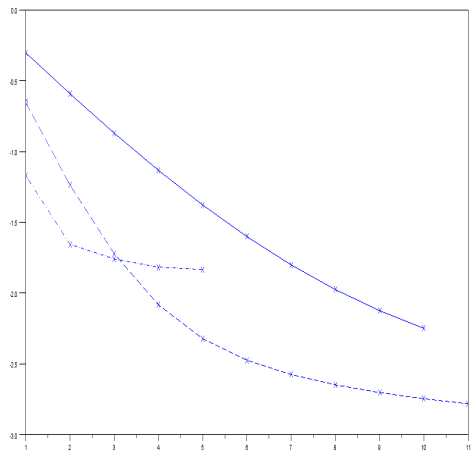}
\caption{\label{half1}Incomplete noiseless data on the half circle from Shepp-Logan phantom with constant speed: relative mean square error plots (BFN in solid line, CG in dashed-dotted line and NS in dashed line).}
\end{figure}

\begin{figure}[!h]\centering
\includegraphics[width=.98\textwidth,height=.98\textwidth]{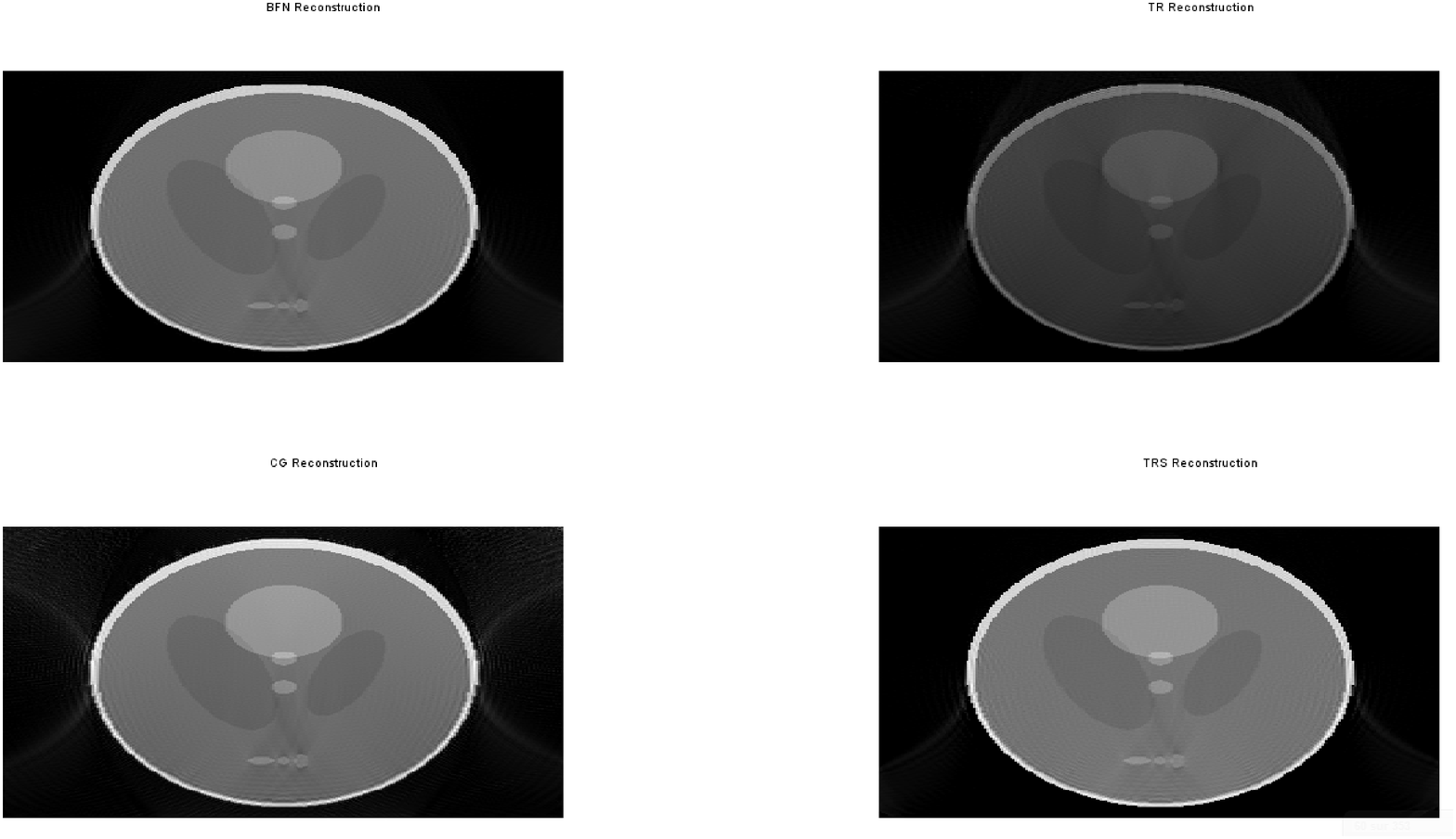}
$$
\begin{array}.{cccc}.
\text{BFN}&\text{CG}&\text{TR}&\text{NS}\\
10.6 (10)&16 (6)&52&6.2 (10)
\end{array}
$$
\caption{\label{half2}Incomplete noiseless data on the half circle from Shepp-Logan phantom with constant speed: reconstructions with BFN (top left), TR (top right), CG (bottom left) and NS (bottom right) and relative mean square errors table (error, number of iterations, limited to 10).}
\end{figure}

In this situation (constant speed and noiseless data on the half circle), the three methods are compared in Fig.~\ref{half1} and~\ref{half2}. NS is excellent, while BFN equivalent to CG convergence in about~5 iterations, but keeps on improving the reconstruction along next iterations. As always, even if its reconstruction is the worst one, CG converges in less iterations and yields a good rough estimate faster than the others.
%\emph{Adding noise}: errors grows highly for BFN, CG and NS. \emph{Attenuating scheme}: very interesting for BFN, NS and CG, less for TR.

%\begin{table}[!h]\centering
%$$
%\begin{array}.{ccc}.
%\text{BFN}&\text{CG}&\text{NS}\\
%10.6 (10)&16 (6)&6.2 (10)
%%\text{$15\%$ noise}&36.4 (8)&55.7 &40 (10)&38.7 (3)\\
%%\e\neq0&23.1 (10)&54.5&26.7 (5)&22.3 (6)
%\end{array}
%$$ 
%\caption{\label{half}Incomplete data, on the half circle.}
%\end{table}

Still testing the robustness of the methods, we reduce the number of sensors, which are still homogeneously located on the observation circle, but only one out of~$\sigma$ is considered in the data. The results shown in Table~\ref{dc} highlight some redundancy aspect of the data since the reconstructions have a similar quality in all situations, particularly with CG which is very robust here. That is obvious at convergence of the algorithms, but one must consider more than~100 iterations for BFN and NS when~$\sigma=8$. Finally, if the intensity is not well recovered, forms are reconstructed with few artifacts.

\begin{table}[!h]\centering
$$
\begin{array}.{lccc}.
\text{Settings}&\text{BFN}&\text{CG}&\text{NS}\\
\sigma=1&2.3 (10)&7.1 (3)&1.4 (10)\\
\sigma=2&5.4 (10)&8.9 (4)&4.5 (10)\\
\sigma=4&12.1 (10)&12.4 (4)&11.5 (10)\\
\sigma=8&25.3 (10)&15 (5)&18.5 (10)
\end{array}
$$ 
\caption{\label{dc}Partial noiseless data (one out of~$\sigma$ sensor records data) from Shepp-Logan phantom with constant speed: relative mean square errors (number of iterations, limited to 10).}
\end{table}

%Final time doubled on Fig.~\ref{time}.
%
%\begin{table}[!h]\centering
%$$
%\begin{array}.{lcccc}.
%\text{Settings}&\text{BFN}&\text{TR}&\text{CG}&\text{NS}\\
%T=\sqrt2&38.8 (2)&42.8&37.7 (6)&35.8 (2)\\
%T=2\sqrt2&38.4 (4)&59.9&37.8 (10)&37.9 (4)
%\end{array}
%$$ 
%\caption{\label{time}Trapping speed TS1 and~$15\%$ noise. Comparison according to the final time.}
%\end{table}

\begin{figure}[!h]\centering
\includegraphics[width=.5\textwidth,height=.5\textwidth]{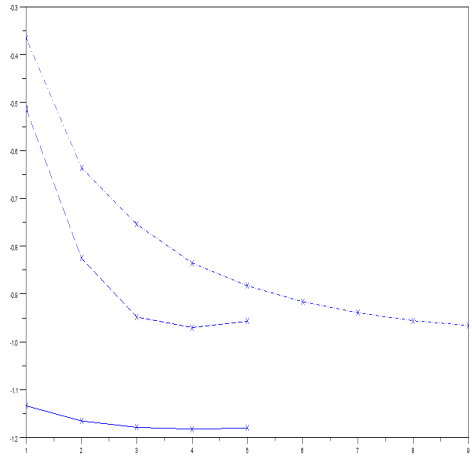}
\caption{\label{bilan1}Incomplete partial data ($\sigma=2$) on the half circle with $15\%$ noise from skull phantom and with constant speed: relative mean square error plots (BFN in solid line, CG in dashed-dotted line and NS in dashed line).}
\end{figure}

\begin{figure}[!h]\centering
\includegraphics[width=.98\textwidth,height=.98\textwidth]{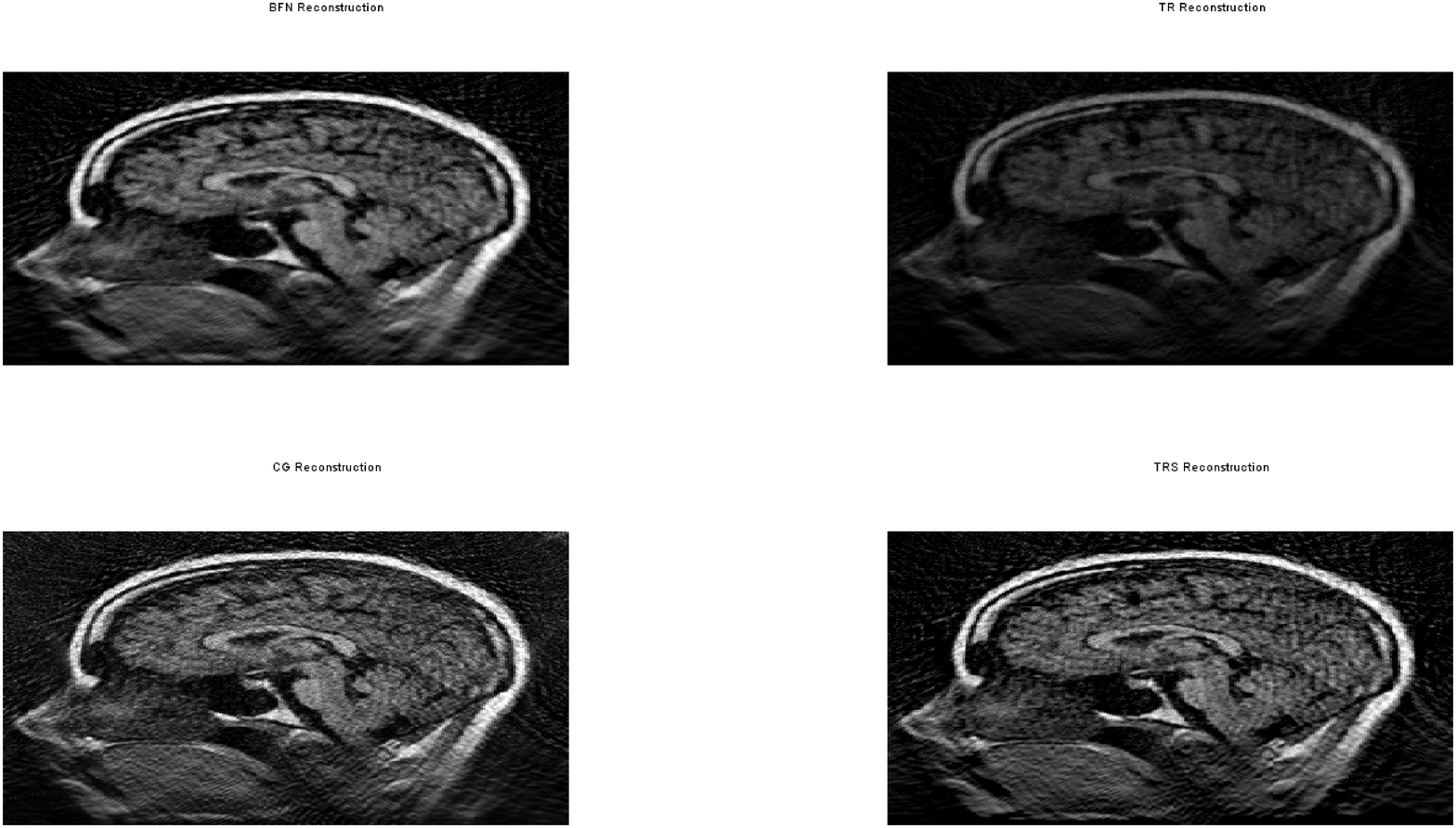}
$$
\begin{array}.{cccc}.
\text{BFN}&\text{CG}&\text{TR}&\text{NS}\\
30.7 (4)&37.8 (10)&59.9&37.9 (4)
\end{array}
$$
\caption{\label{bilan2}Incomplete $15\%$ noise data on the half circle from skull phantom with constant speed and partial data ($\sigma=2$): reconstructions with BFN (top left), TR (top right), CG (bottom left) and NS (bottom right) and relative mean square errors table (error, number of iterations, limited to 10).}
\end{figure}

Then we sum up the different unfavorable conditions from above and show the results on Fig.~\ref{bilan1} and~\ref{bilan2}. Here data are acquired on the half circle, only one out of two sensors is considered, a ~$15\%$ level noise is added and the speed is constant. The last object, a skull, has been chosen in this last experiment in order to consider a more complex structured object.

Finally, even if the reconstructions obtained in this last situation are not satisfying, it may fit with very bad real conditions (where the noise level can be far greater than~$15\%$ for example) and allows us to test the robustness of the methods. In this view, Table~\ref{vvarpibr} shows the errors if we consider variable speeds in a similar case but with~$\sigma=1$.
%Noisy incomplete half circle data on Table~\ref{noiseinc}.

\begin{table}[!h]\centering
$$
\begin{array}.{lcccc}.
\text{Speed map}&\text{BFN}&\text{CG}&\text{NS}\\
c\equiv1&36.4 (8)&40 (10)&38.7 (3)\\
\text{NTS}&36.3 (8)&40.6 (10)&38.6 (3)\\
\text{TS1}&46.3 (4)&50.9 (8)&51.9 (3)\\
\text{TS2}&39.1 (8)&42.3 (9)&41.5 (4)
\end{array}
$$ 
\caption{\label{vvarpibr}Incomplete $15\%$ noise data on the half circle from Shepp-Logan phantom: relative mean square errors (number of iterations, limited to 10).}
\end{table}

%%%%%%%%%%%%%%%%%%%%%%%%%%%%%%%%%%%%%%%%%%%%%%%%%%%%%%%%%%%%%
\subsection{Implementation conclusions}%%%%%%%%%%%%%%%%%%%%%%
%%%%%%%%%%%%%%%%%%%%%%%%%%%%%%%%%%%%%%%%%%%%%%%%%%%%%%%%%%%%%

CG reacts much badly to noise addition above all, but it provides good reconstructions after a couple of iterations in comparison with BFN and NS, so that it can be used to get a rough estimate for another method. Moreover CG keeps interesting when data are incomplete, not overreacting. On the whole, BFN and NS provide good and quite equivalent reconstructions with a similar cost. Finally, NS furnishes better reconstruction, but appears less robust than BFN to noise and to restriction of the number of sensors. Numerical attenuation provides a good regularization that offsets noise for NS and BFN. Notice that this artificial numerical attenuation may allow us to consider lossy medium in back and forth evolutions as the physical loss does not exceed the artificial one.

Depending on whether data are complete or not, one may choose CG to get a first rough estimate, then NS or BFN to reconstruct the object. It has to be said that not stopping iterations to~10 often improves results in noiseless and/or complete data situations, but that this number of iterations provides the best estimate in less favorable settings.

%%%%%%%%%%%%%%%%%%%%%%%%%%%%%%%%%%%%%%%%%%%%%%%%%%%%%%%%%%%%%
\section{Conclusions}%%%%%%%%%%%%%%%%%%%%%%%%%%%%%%%%%%%%%%
%%%%%%%%%%%%%%%%%%%%%%%%%%%%%%%%%%%%%%%%%%%%%%%%%%%%%%%%%%%%%%
We have established the theoretical geometrical convergence of the TAT-BFN algorithm, with an explicit convergence rate in an ideal situation. In principle, this method enables to make weak assumptions in the TAT modeling (variable sound speed, damping, incomplete data, etc.), since the model is considered as a weak constraint. Even though the proof of the convergence has only been stated in an ideal framework, the numerical implementation of the algorithm showed that the BFN method remains efficient in more general situations. Indeed this technique still yields robust reconstructions in the noisy and incomplete data situation. Even though the convergence rates of the Conjugate Gradient and of the Neumann Series methods appear to be better, the BFN provides more satisfactory asymptotic behavior and minimal relative error, along with a low numerical complexity, in quite unfavorable observation situations. 

From these preliminary results, one should now focus on their generalization to more realistic wave equations. For example, one may introduce some damping
in the model and study the theoretical validity and convergence of the methods. Our first attempts in this direction (mainly numerical implementations) 
suggest future interesting developments.

%\nocite{*}

\printbibliography

\end{document}